\documentclass[12pt]{article}
\usepackage{amssymb,amsmath,amsthm,amsfonts,enumerate, bbm, mathdots}

\usepackage{latexsym}
\usepackage{amscd}
\usepackage{stmaryrd}
\usepackage{color}
\usepackage[all]{xypic}
\usepackage{epsfig}
\usepackage{graphics}
\usepackage{ifthen}
\usepackage{varioref}
\usepackage{rotating}
\usepackage{extarrows}
\usepackage{cite}

\numberwithin{equation}{section}

\theoremstyle{plain}
\newtheorem{thm}{Theorem}[section]
\newtheorem{cor}[thm]{Corollary}
\newtheorem{lem}[thm]{Lemma}
\newtheorem{prop}[thm]{Proposition}

\definecolor{darkgreen}{rgb}{0.0625,0.64,0.0625}
\usepackage[
            pdfstartview=FitH,
            CJKbookmarks=true,
            bookmarksnumbered=true,
            bookmarksopen=true,
            colorlinks,
            pdfborder=001,
            linkcolor=blue,
            anchorcolor=green,
            citecolor=red
            ]{hyperref}

\newfont{\scyr}{wncyr10 scaled 550}
\def\shuffle{\,\mbox{\bf \scyr X}\,}
\def\tshuffle{\overset{t}{\shuffle}}
\def\tast{\overset{t}{\ast}}

\def\proof{\noindent {\bf Proof.\;}}

\def\wt{\operatorname{wt}}
\def\dep{\operatorname{dep}}
\def\height{\operatorname{ht}}
\def\reg{\operatorname{reg}}

\allowdisplaybreaks

\begin{document}

\title{Algebraic relations of interpolated multiple zeta values}

\date{\small ~ \qquad\qquad School of Mathematical Sciences, Tongji University \newline No. 1239 Siping Road,
Shanghai 200092, China}

\author{Zhonghua Li \thanks{E-mail address: zhonghua\_li@tongji.edu.cn}~\thanks{The author is supported by the National Natural Science Foundation of
China (Grant No. 11471245).}
}

\maketitle

\begin{abstract}
Interpolated multiple zeta values can be regarded as interpolation polynomials of multiple zeta values and multiple zeta-star values. In this paper, we give some algebraic relations of interpolated multiple zeta values, such as the symmetric sum formula, the shuffle regularized sum formula, a weighted sum formula and some evaluation formulas with even arguments. All the algebraic relations provided in this paper are deduced from the extended double shuffle relations.
\end{abstract}

{\small

}


\section{Introduction}\label{Sec:Intro}

To study the multiple zeta values and the multiple zeta-star values simultaneously, S. Yamamoto introduced an interpolation polynomial of these two types of values in \cite{Yamamoto}. We call a finite sequence of positive integers an index.  An index $\mathbf{k}=(k_1,k_2,\ldots,k_n)$ is admissible if $k_1\geqslant 2$. For an admissible index $\mathbf{k}=(k_1,k_2,\ldots,k_n)$ and for a variable $t$, the interpolated multiple zeta value ($t$-MZV) $\zeta^t(\mathbf{k})$ is defined as
$$\zeta^t(\mathbf{k})=\zeta^t(k_1,\ldots,k_n)=\sum\limits_{\mathbf{p}=(k_1\Box k_2\Box \cdots\Box k_n)\atop \Box=\lq\lq," \text{\, or\,}\lq\lq+"}t^{n-\dep(\mathbf{p})}\zeta(\mathbf{p})(\in\mathbb{R}[t]).$$
Here for any index $\mathbf{k}=(k_1,k_2,\ldots,k_n)$, the depth of $\mathbf{k}$ is defined by $\dep(\mathbf{k})=n$, and if $\mathbf{k}$ is admissible, then the multiple zeta value (MZV) $\zeta(\mathbf{k})$ is defined by
$$\zeta(\mathbf{k})=\zeta(k_1,\ldots,k_n)=\sum\limits_{m_1>m_2>\cdots>m_n>0}\frac{1}{m_1^{k_1}m_2^{k_2}\cdots m_n^{k_n}}.$$
It is obviously that $\zeta^0(\mathbf{k})=\zeta(\mathbf{k})$. And one can show that
$$\zeta^1(\mathbf{k})=\zeta^{\star}(\mathbf{k})=\zeta^{\star}(k_1,\ldots,k_n)=\sum\limits_{m_1\geqslant m_2\geqslant\cdots\geqslant m_n>0}\frac{1}{m_1^{k_1}m_2^{k_2}\cdots m_n^{k_n}},$$
which is called a multiple zeta-star value (MZSV).

Many relations, which are satisfied by MZVs and MZSVs simultaneously, are generalized to $t$-MZVs. A typical example is the sum formula which was first established by S. Yamamoto in \cite{Yamamoto}. For an index $\mathbf{k}=(k_1,k_2,\ldots,k_n)$, besides the depth, we define the weight of $\mathbf{k}$ by
$$\wt(\mathbf{k})=k_1+k_2+\cdots+k_n.$$
Denote by $I_{k,n}$ the set of all indices with weight $k$ and depth $n$, and by $I^0_{k,n}$ the subset of $I_{k,n}$ containing all admissible indices. Then the sum formula of $t$-MZVs claims that
\begin{align}
\sum\limits_{\mathbf{k}\in I^0_{k,n}}\zeta^t(\mathbf{k})=\left(\sum\limits_{i=0}^{n-1}\binom{k-1}{i}t^i(1-t)^{n-1-i}\right)\zeta(k),
\label{Eq:SumFormula}
\end{align}
where $k$ and $n$ are positive integers with $k>n$. We list some other known relations of $t$-MZVs
\begin{itemize}
  \item the cyclic sum formula proved by S. Yamamoto in \cite{Yamamoto};
  \item Kawashima's relations obtained by T. Tanaka and N. Wakabayashi in \cite{Tanaka-Wakabayashi};
  \item the extended double shuffle relations and Hoffman's relation established by N. Wakabayashi in \cite{Wakabayashi} and by C. Qin and the author in \cite{Li-Qin2017} independently.
\end{itemize}

In this paper, we focus on the algebraic aspects  of $t$-MZVs, especially the algebraic relations which are deduced from the extended double shuffle relations. For this purpose, we first recall the algebraic setting of the extended double shuffle relations in Section \ref{Sec:EDS}. Different from \cite{Li-Qin2017} and \cite{Wakabayashi}, we deal with a more general situation and obtain some equivalent statements. In fact, we find that each equivalent statement of the extended double shuffle relations proposed in \cite{Li-Qin} either in terms of MZVs or in terms of MZSVs has a $t$-MZVs version. Then we give some algebraic relations of $t$-MZVs in Section \ref{Sec:AlgRel}. We provide some new relations of $t$-MZVs, such as the symmetric sum formula, the shuffle regularized sum formula and a weighted sum formula which generalizes the weighted sum formula of L. Guo and B. Xie given in \cite{Guo-Xie2009}. Finally, in Section \ref{Sec:EvaFormula} we give some evaluation formulas of $t$-MZVs, including the formula of $\zeta^t(2k,\ldots,2k)$ and a restricted sum formula which evaluates the sum
$$\sum\limits_{k_1+\cdots+k_n=k\atop k_i\geqslant 1}\zeta^t(2mk_1,\ldots,2mk_n).$$


\section{Extended double shuffle relations}\label{Sec:EDS}

In this section, we recall the algebraic setting of the extended double shuffle relations of $t$-MZVs from \cite{Li-Qin2017,Wakabayashi}. While different from \cite{Li-Qin2017,Wakabayashi}, we find that the extended double shuffle relations of $t$-MZVs are in fact equivalent to the extended double shuffle relations of MZVs. Therefore as in \cite{Li-Qin}, we obtain some equivalent statements systematically and deal with a more general situation.

Let $\mathbb{K}$ be a field of characteristic zero and $R$ be a commutative $\mathbb{K}$-algebra with unitary.

\subsection{Algebraic setup}

The algebraic setup of $t$-MZVs was first established in \cite{Tanaka-Wakabayashi,Yamamoto}. In the case of $t=0$, one obtains the classical algebraic setup of MZVs as in \cite{Hoffman1997,Hoffman-Ohno}, and in the case of $t=1$, one gets the algebraic setup of MZSVs as in \cite{Muneta}.

Let $A=\{x,y\}$ be an alphabet which contains two noncommutative letters, and let $A^{\ast}$ be the set of all words generated by $A$. We denote by $1$ the empty word. For a word $w\in A^{\ast}$ and a letter $a\in A$, let us denote by $d_a(w)$ the number of $a$'s contained in $w$, and set $|w|=d_x(w)+d_y(w)$. Let $\mathfrak{h}_t=\mathbb{K}[t]\langle A\rangle$ be the noncommutative polynomial algebra over $\mathbb{K}[t]$ generated by $A$. There are two subalgebras of $\mathfrak{h}_t$
$$\mathfrak{h}_t^1=\mathbb{K}[t]+\mathfrak{h}_ty,\quad \mathfrak{h}_t^0=\mathbb{K}[t]+x\mathfrak{h}_ty.$$

The $t$-shuffle product $\tshuffle$ on $\mathfrak{h}_t$ is $\mathbb{K}[t]$-bilinear, and satisfies the rules
\begin{itemize}
  \item[(S1)] $1\tshuffle w=w\tshuffle 1=w$,
  \item[(S2)] $aw_1\tshuffle bw_2=a(w_1\tshuffle bw_2)+b(aw_1\tshuffle w_2)-\delta(w_1)\rho(a)bw_2-\delta(w_2)\rho(b)aw_1$,
\end{itemize}
where $w,w_1, w_2\in A^{\ast}$, $a, b\in A$, the map $\delta: A^{\ast}\rightarrow\{0, 1\}$ is defined by
$$\delta(w)=\begin{cases}
1 & \text{\;if\;} w=1,\\
0 & \text{\;if\;} w\neq 1,
\end{cases}$$
and the map $\rho: A\rightarrow \mathfrak{h}_t$ is defined by
$$\rho(x)=0,\quad  \rho(y)=tx.$$
Then one can prove that the $t$-shuffle product $\tshuffle$ is commutative and associative. Under this new product, $\mathfrak{h}_t$ becomes a commutative and associative algebra, $\mathfrak{h}_t^1$ and $\mathfrak{h}_t^0$ are still subalgebras of $\mathfrak{h}_t$.

For any $k\in\mathbb{N}$, where $\mathbb{N}$ is the set of positive integers, set $z_k=x^{k-1}y$. The $t$-harmonic shuffle (stuffle) product on $\mathfrak{h}_t^1$ is  $\mathbb{K}[t]$-bilinear, and satisfies the rules
\begin{itemize}
  \item[(H1)] $1\tast w=w\tast 1=w$,
  \item[(H2)] $z_kw_1\tast z_lw_2=z_k(w_1\tast z_lw_2)+z_l(z_kw_1\tast w_2)+(1-2t)z_{k+l}(w_1\tast w_2)+[1-\delta(w_1)\delta(w_2)](t^2-t)x^{k+l}(w_1\tast w_2)$,
\end{itemize}
where $w,w_1,w_2\in A^{\ast}\cap\mathfrak{h}_t^1$ and $k,l\in\mathbb{N}$. As for the $t$-shuffle product, one can prove that the $t$-harmonic shuffle product $\tast$ is commutative and associative. Hence under this new product, $\mathfrak{h}_t^1$ becomes a commutative and associative algebra, and $\mathfrak{h}_t^0$ is still a subalgebra.

In the case of $t=0$, we simply denote $\mathfrak{h}_0$, $\mathfrak{h}_0^1$ and $\mathfrak{h}_0^0$ by $\mathfrak{h}$, $\mathfrak{h}^1$ and $\mathfrak{h}^0$, respectively. Note that we have
$$\mathfrak{h}_t=\mathfrak{h}[t],\quad \mathfrak{h}_t^1=\mathfrak{h}^1[t],\quad \mathfrak{h}_t^0=\mathfrak{h}^0[t].$$
We also denote $\overset{0}{\shuffle}$ by $\shuffle$, which is the usual shuffle product. And one can extend $\shuffle$ to the whole space $\mathfrak{h}_t$, such that $\shuffle$ is $\mathbb{K}[t]$-bilinear. Similarly, $\ast=\overset{0}{\ast}$ is the usual harmonic shuffle product, and this product can be extended to $\mathfrak{h}_t^1$ by $\mathbb{K}[t]$-bilinearities.

Under the algebraic setup, the most important fact is that one can associate the general case with the special case of $t=0$. Let $\sigma_t$ be an automorphism of the noncommutative algebra $\mathfrak{h}_t$ determined by
$$\sigma_t(x)=x,\quad \sigma_t(y)=tx+y.$$
Note that $\sigma_t^{-1}=\sigma_{-t}$. The $\mathbb{K}[t]$-linear map $S_t:\mathfrak{h}_t\rightarrow\mathfrak{h}_t$ is defined by $S_t(1)=1$ and
$$S_t(wa)=\sigma_t(w)a,\quad (\forall w\in A^{\ast}, \forall a\in A).$$
Note that $S_t$ is invertible, and $S_t^{-1}=S_{-t}$. Moreover, we have $S_t(\mathfrak{h}_t^1)=\mathfrak{h}_t^1$ and $S_t(\mathfrak{h}_t^0)=\mathfrak{h}_t^0$. In the absence of confusion, both $S_t|_{\mathfrak{h}_t^1}$ and $S_t|_{\mathfrak{h}_t^0}$ are simply denoted by $S_t$. Then for any $w_1,w_2\in\mathfrak{h}_t$, we have
$$S_t(w_1\tshuffle w_2)=S_t(w_1)\shuffle S_t(w_2),\quad S_{-t}(w_1\shuffle w_2)=S_{-t}(w_1)\tshuffle S_{-t}(w_2).$$
And for any $w_1,w_2\in\mathfrak{h}_t^1$, we have
$$S_t(w_1\tast w_2)=S_t(w_1)\ast S_t(w_2),\quad S_{-t}(w_1\ast w_2)=S_{-t}(w_1)\tast S_{-t}(w_2).$$

For the later use, we prepare a lemma below.

\begin{lem}\label{Lem:S-Com}
For variables $s$ and $t$, we have
$$\sigma_s\circ \sigma_t=\sigma_{s+t},\quad S_s\circ S_t=S_{s+t}.$$
\end{lem}

\proof
Since
$$(\sigma_s\circ \sigma_t)(x)=\sigma_s(x)=x=\sigma_{s+t}(x)$$
and
$$(\sigma_s\circ \sigma_t)(y)=\sigma_s(tx+y)=tx+sx+y=\sigma_{s+t}(y),$$
we get $\sigma_s\circ \sigma_t=\sigma_{s+t}$. Then for any $a\in A$ and $w\in A^{\ast}$, we have
$$(S_s\circ S_t)(wa)=S_s(\sigma_t(w)a)=\sigma_s(\sigma_t(w))a=\sigma_{s+t}(w)a=S_{s+t}(wa),$$
which induces that $S_s\circ S_t=S_{s+t}$.
\qed


\subsection{Extended double shuffle relations}

Let $Z_R:\mathfrak{h}^0\rightarrow R$ be a $\mathbb{K}$-linear map, which can be extended to a $\mathbb{K}[t]$-linear map $Z_R:\mathfrak{h}_t^0=\mathfrak{h}^0[t]\rightarrow R[t]$. Define
$$Z_R^t=Z_R\circ S_t:\mathfrak{h}_t^0\rightarrow R[t],$$
then $Z_R^t$ is $\mathbb{K}[t]$-linear. Now for any admissible index $\mathbf{k}=(k_1,\ldots,k_n)$, we define the MZV and the $t$-MZV associated with $Z_R$ respectively by
$$\zeta_R(\mathbf{k})=Z_R(z_{k_1}\cdots z_{k_n}),\quad \zeta_R^t(\mathbf{k})=Z_R^{t}(z_{k_1}\cdots z_{k_n}).$$

As in \cite{Ihara-Kaneko-Zagier,Li-Qin}, if the $\mathbb{K}$-linear map $Z_R: \mathfrak{h}^0\rightarrow R$ satisfies the conditions
$$Z_R(w_1\shuffle w_2)=Z_R(w_1)Z_R(w_2)=Z_R(w_1\ast w_2), \quad (\forall w_1,w_2\in\mathfrak{h}^0),$$
we call the map $Z_R$ satisfying the (finite) double shuffle relations. It is easy to show that

\begin{prop}\label{Prop:FDS}
For any $\mathbb{K}$-linear map $Z_R:\mathfrak{h}^0\rightarrow R$, $Z_R$ satisfies the double shuffle relations if and only if
$$Z_R^t(w_1\tshuffle w_2)=Z_R^t(w_1)Z_R^t(w_2)=Z_R^t(w_1\tast w_2), \quad (\forall w_1,w_2\in\mathfrak{h}^0_t).$$
\end{prop}

Below, if $Z_R:\mathfrak{h}^0\rightarrow R$ satisfies the double shuffle relations, we always assume that $Z_R(1)=1$. Hence if $\mathbf{k}$ is the empty index, we have $\zeta_R(\mathbf{k})=\zeta_R^t(\mathbf{k})=1$.

Now we recall the extended double shuffle relations. Under the shuffle product, $(\mathfrak{h}^1,\shuffle)=(\mathfrak{h}^0,\shuffle)[y]$ is a polynomial algebra. More precisely, for any $w\in \mathfrak{h}^1$, there exist $w_i\in\mathfrak{h}^0$ such that
$$w=\sum\limits_{i\geqslant 0}w_i\shuffle y^{\shuffle i}.$$
Hence there is a map $\reg_{\shuffle}:\mathfrak{h}^1\rightarrow\mathfrak{h}^0$ defined by $\reg_{\shuffle}(w)=w_0$. The map $\reg_{\shuffle}:(\mathfrak{h}^1,\shuffle)\rightarrow(\mathfrak{h}^0,\shuffle)$ is an algebra homomorphism and is called the shuffle regularization. Similarly, since $(\mathfrak{h}^1,\ast)=(\mathfrak{h}^0,\ast)[y]$, one can define the harmonic shuffle regularization map $\reg_{\ast}:(\mathfrak{h}^1,\ast)\rightarrow(\mathfrak{h}^0,\ast)$, which is also an algebra homomorphism.

Let $Z_R:\mathfrak{h}^0\rightarrow R$ be a $\mathbb{K}$-linear map which satisfies the double shuffle relations. Then there exist unique algebra homomorphisms
$$Z_R^{\shuffle}:(\mathfrak{h}^1,\shuffle)\rightarrow R[T],\quad Z_R^{\ast}:(\mathfrak{h}^1,\ast)\rightarrow R[T],$$
which are determined by
$$Z_R^{\shuffle}|_{\mathfrak{h}^0}=Z_R=Z_R^{\ast}|_{\mathfrak{h}^0},\quad Z_R^{\shuffle}(y)=T=Z_R^{\ast}(y),$$
where $T$ is a variable. Then the following items are equivalent which is proved in \cite{Ihara-Kaneko-Zagier}:
\begin{itemize}
    \item[(i)] $(Z_{R}^{\shuffle}-\rho_{R}\circ Z_{R}^{\ast})(w_1)=0$ for all $w_1\in\mathfrak{h}^{1}$;
    \item[(ii)] $(Z_{R}^{\shuffle}-\rho_{R}\circ Z_{R}^{\ast})(w_1)|_{T=0}=0$ for all $w_1\in\mathfrak{h}^{1}$;
    \item[(iii)] $Z_{R}^{\shuffle}(w_1\shuffle w_0-w_1\ast w_0)=0$ for all $w_1\in\mathfrak{h}^{1}$ and all $w_0\in\mathfrak{h}^{0}$;
    \item[(iii$'$)] $Z_{R}^{\ast}(w_1\shuffle w_0-w_1\ast w_0)=0$ for all $w_1\in\mathfrak{h}^{1}$ and all $w_0\in\mathfrak{h}^{0}$;
    \item[(iv)] $Z_{R}(\reg_{\shuffle}(w_1\shuffle w_0-w_1\ast w_0))=0$ for all $w_1\in\mathfrak{h}^{1}$ and all $w_0\in\mathfrak{h}^{0}$;
    \item[(iv$'$)] $Z_{R}(\reg_{\ast}(w_1\shuffle w_0-w_1\ast w_0))=0$ for all $w_1\in\mathfrak{h}^{1}$ and all $w_0\in\mathfrak{h}^{0}$;
    \item[(v)] $Z_{R}(\reg_{\shuffle}(y^m\ast w_0))=0$ for all $m\in\mathbb{N}$ and all $w_0\in \mathfrak{h}^{0}$;
    \item[(v$'$)] $Z_{R}(\reg_{\ast}(y^m\shuffle w_0-y^m\ast w_0))=0$ for all $m\in\mathbb{N}$ and all $w_0\in \mathfrak{h}^{0}$;
    \item[(vi)] $Z_{R}(\partial_n(w_0))=0$ for all $n\in\mathbb{N}$ and all $w_0\in \mathfrak{h}^{0}$;
    \item[(vii)] $Z_{R}((\sigma_{m}-\overline{\sigma_{m}})(w_0))=0$ for all $m\in\mathbb{Z}_{\geqslant 0}$ and all $w_0\in \mathfrak{h}^{0}$.
\end{itemize}
We call a $\mathbb{K}$-linear map $Z_R:\mathfrak{h}^0\rightarrow R$ satisfying the extended double shuffle relations if it satisfies the double shuffle relations and the above equivalent properties. Here $\rho_{R}:R[T]\rightarrow R[T]$ is an $R$-module homomorphism defined by
$$\rho_R(e^{Tu})=A_{R}(u)e^{Tu}$$
with $u$ a variable and
$$A_{R}(u)=\exp\left(\sum\limits_{n=2}^{\infty}\frac{(-1)^{n}}{n}\zeta_R(k)u^n\right)\in R[[u]].$$
The map $\partial_n$ is the derivation on $\mathfrak{h}$ determined by
$$\partial_n(x)=x(x+y)^{n-1}y,\quad \partial_n(y)=-x(x+y)^{n-1}y.$$
The map $\sigma_m:\mathfrak{h}^1\rightarrow\mathfrak{h}^1$ is a $\mathbb{K}$-linear map defined by $\sigma_m(1)=1$ and
$$\sigma_m(z_{k_1}\cdots z_{k_n})=\sum\limits_{\varepsilon_1+\cdots+\varepsilon_n=m\atop \varepsilon_i\geqslant 0}z_{k_1+\varepsilon_1}\cdots z_{k_n+\varepsilon_n},\quad (\forall n,k_1,\ldots,k_n\in\mathbb{N}).$$
Finally, $\overline{\sigma_{m}}=\tau\circ\sigma_m\circ\tau$, where $\tau$ is the antiautomorphism of the noncommutative algebra $\mathfrak{h}$ determined by
$$\tau(x)=y,\qquad \tau(y)=x,$$
and $\mathbb{Z}_{\geqslant 0}$ is the set of nonnegative integers. We also recall that a $\mathbb{K}$-linear map $\mathcal{D}:\mathfrak{h}\rightarrow \mathfrak{h}$ is a derivation if it satisfies
$$\mathcal{D}(w_1w_2)=\mathcal{D}(w_1)w_2+w_1\mathcal{D}(w_2),\quad (\forall w_1,w_2\in\mathfrak{h}).$$

Using the map $S_t$, every equivalent properties above can be stated in terms of $t$-MZVs. As shown in \cite{Li-Qin2017,Wakabayashi}, as $t$-shuffle algebras, we have $\mathfrak{h}_t^1=\mathfrak{h}_t^0[y]$. Hence we have an algebra homomorphism $\reg_{\tshuffle}:(\mathfrak{h}_t^1,\tshuffle)\rightarrow(\mathfrak{h}_t^0,\tshuffle)$. Similarly, as $t$-harmonic shuffle algebras, we have $\mathfrak{h}_t^1=\mathfrak{h}_t^0[y]$, which induces an algebra homomorphism $\reg_{\tast}:(\mathfrak{h}_t^1,\tast)\rightarrow(\mathfrak{h}_t^0,\tast)$. It is easy to prove that
$$\reg_{\tshuffle}=S_{-t}\circ \reg_{\shuffle}\circ S_t,\quad \reg_{\tast}=S_{-t}\circ \reg_{\ast}\circ S_t,$$
where $\reg_{\shuffle}$ and $\reg_{\ast}$ are extended to $\mathbb{K}[t]$-linear maps from $\mathfrak{h}_t^1=\mathfrak{h}^1[t]$ to $\mathfrak{h}_t^0=\mathfrak{h}^0[t]$. Furthermore, for a $\mathbb{K}$-linear map $Z_R:\mathfrak{h}^0\rightarrow R$ which satisfies the double shuffle relations, there exist unique algebra homomorphisms $Z_R^{t,\tshuffle}:(\mathfrak{h}^1_t,\tshuffle)\rightarrow R[t,T]$ and $Z_R^{t,\tast}:(\mathfrak{h}^1_t,\tast)\rightarrow R[t,T]$ such that
$$Z_R^{t,\tshuffle}|_{\mathfrak{h}_t^0}=Z_R^t=Z_R^{t,\tast}|_{\mathfrak{h}_t^0},\quad Z_R^{t,\tshuffle}(y)=T=Z^{t,\tast}_R(y).$$
Extending the maps $Z_R^{\shuffle}$ and $Z_R^{\ast}$ to be $\mathbb{K}[t]$-linear maps from $\mathfrak{h}_t^1$ to $R[t,T]$, we get
$$Z_R^{t,\tshuffle}=Z_R^{\shuffle}\circ S_t,\quad Z_R^{t,\tast}=Z_R^{\ast}\circ S_t.$$

For any $n\in\mathbb{N}$, the derivation $\partial_n$ on $\mathfrak{h}$ can be extended to a $\mathbb{K}[t]$-linear map on $\mathfrak{h}_t$. Let us define $\partial_n^t=S_{-t}\circ \partial_n\circ S_t$, then as in \cite{Li-Qin}, we can show that $\partial_n^t$ is a left $\widetilde{S_t}$-derivation on $\mathfrak{h}_t$ and satisfies
$$\partial_n^t(x)=x(x-tx+y)^{n-1}y,\quad \partial_n^t(y)=-x(x-tx+y)^{n-1}y.$$
Here $\widetilde{S_t}=S_{-t}\circ \sigma_t$ is $\mathbb{K}[t]$-linear, which satisfies $\widetilde{S_t}(1)=1$ and
$$\widetilde{S_t}(wa)=w\sigma_t(a),\quad (\forall w\in\mathfrak{h}_t, \forall a\in A).$$
We recall from \cite{Li-Qin}, for $\mathbb{K}[t]$-linear maps $\mathcal{D}$ and $\mathcal{L}$ on $\mathfrak{h}_t$ with $\mathcal{D}(1)=0$ and $\mathcal{L}$ invertible, if
$$\mathcal{D}(w_1w_2)=(\mathcal{L}^{-1}\circ \mathcal{D}\circ\mathcal{L})(w_1)w_2+w_1\mathcal{D}(w_2),\quad (1\neq \forall w_1,w_2\in A^{\ast}),$$
we call $\mathcal{D}$ is a left $\mathcal{L}$-derivation on $\mathfrak{h}_t$.

Finally, for any $m\in\mathbb{Z}_{\geqslant 0}$, we extend $\sigma_m$ and $\overline{\sigma_m}$ to be $\mathbb{K}[t]$-linear maps on $\mathfrak{h}_t^1$, and define
$$\sigma_m^t=S_{-t}\circ \sigma_m\circ S_t,\quad \overline{\sigma_m}^t=S_{-t}\circ \overline{\sigma_m}\circ S_t.$$
Then after extending $\rho_R$ to be a $R[t]$-module endomorphism on $R[t,T]$, as in \cite{Li-Qin}, we can prove the following theorem.

\begin{thm}\label{Thm:EDS-tMZV}
Assume that the $\mathbb{K}$-linear map $Z_R:\mathfrak{h}^0\rightarrow R$ satisfies the double shuffle relations. Then the following properties are equivalent:
\begin{description}
    \item[(0)] $Z_R$ satisfies the extended double shuffle relations;
    \item[(1)] $(Z_{R}^{\tshuffle}-\rho_{R}\circ Z_{R}^{\tast})(w_1)=0$ for all $w_1\in\mathfrak{h}^{1}_t$;
    \item[(2)] $(Z_{R}^{\tshuffle}-\rho_{R}\circ Z_{R}^{\tast})(w_1)|_{T=0}=0$ for all $w_1\in\mathfrak{h}^{1}_t$;
    \item[(3)] $Z_{R}^{\tshuffle}(w_1\tshuffle w_0-w_1\tast w_0)=0$ for all $w_1\in\mathfrak{h}^{1}_t$ and all $w_0\in\mathfrak{h}^{0}_t$;
    \item[(3$'$)] $Z_{R}^{\tast}(w_1\tshuffle w_0-w_1\tast w_0)=0$ for all $w_1\in\mathfrak{h}^{1}_t$ and all $w_0\in\mathfrak{h}^{0}_t$;
    \item[(4)] $Z_{R}^{t}(\reg_{\tshuffle}(w_1\tshuffle w_0-w_1\tast w_0))=0$ for all $w_1\in\mathfrak{h}^{1}_t$ and all $w_0\in\mathfrak{h}^{0}_t$;
    \item[(4$'$)] $Z_{R}^{t}(\reg_{\tast}(w_1\tshuffle w_0-w_1\tast w_0))=0$ for all $w_1\in\mathfrak{h}^{1}_t$ and all $w_0\in\mathfrak{h}^{0}_t$;
    \item[(5)] $Z_{R}^{t}(\reg_{\tshuffle}(y^m\tshuffle w_0-y^m\tast w_0))=0$ for all $m\in\mathbb{N}$ and all $w_0\in \mathfrak{h}^{0}_t$;
    \item[(5$'$)] $Z_{R}^{t}(\reg_{\tast}(y^m\tshuffle w_0-y^m\tast w_0))=0$ for all $m\in\mathbb{N}$ and all $w_0\in \mathfrak{h}^{0}_t$;
    \item[(6)] $Z_{R}^{t}(\partial_n^{t}(w_0))=0$ for all $n\in\mathbb{N}$ and all $w_0\in \mathfrak{h}^{0}_t$;
    \item[(7)] $Z_{R}^{t}((\sigma_{m}^{t}-\overline{\sigma_{m}}^{t})(w_0))=0$ for all $m\in\mathbb{Z}_{\geqslant 0}$ and all $w_0\in \mathfrak{h}^{0}_t$.
\end{description}
\end{thm}

Let $\mathbb{K}=\mathbb{Q}$, $R=\mathbb{R}$, $Z_R=Z:\mathfrak{h}^0\rightarrow \mathbb{R}$ defined by $Z(1)=1$ and
$$Z(z_{k_1}\cdots z_{k_n})=\zeta(k_1,\ldots,k_n),\quad (n,k_1,\ldots,k_n\in\mathbb{N}, k_1\geqslant 2),$$
then the map $Z$ satisfies the extended double shuffle relations by \cite{Ihara-Kaneko-Zagier}. Obviously, $Z^t=Z\circ S_t$ satisfies
$$Z^t(z_{k_1}\cdots z_{k_n})=\zeta^t(k_1,\ldots,k_n),\quad (n,k_1,\ldots,k_n\in\mathbb{N}, k_1\geqslant 2).$$
Hence Proposition \ref{Prop:FDS} and Theorem \ref{Thm:EDS-tMZV} express the extended double shuffle relations in terms of $t$-MZVs. In other words, we get the extended double shuffle relations of $t$-MZVs obtained in \cite{Li-Qin2017,Wakabayashi} and give more equivalent statements.


\section{Some algebraic relations}\label{Sec:AlgRel}

In this section, we give some algebraic relations of $t$-MZVs, which are deduced from the extended double shuffle relations. We first work out the wanted identity in $\mathfrak{h}_t$, then we can obtain the corresponding algebraic relations of $t$-MZVs after applying the map $Z$ or $Z^t$ and the extended double shuffle relations. For simplicity, we assume that $\mathbb{K}=\mathbb{Q}$, $R=\mathbb{R}$ and $Z_R=Z$ from now on.

\subsection{Symmetric sum formula}

The symmetric sum formulas of MZVs and of MZSVs were first established by M. E. Hoffman in \cite[Theorem2.1 and Theorem 2.2]{Hoffman1992}. Here we generalize them to $t$-MZVs. First we introduce some notations. For any $m\in\mathbb{N}$, set
$$c_m(t)=(m-1)!\left[t^m-(t-1)^m\right]=(m-1)!\sum\limits_{i=0}^{m-1}(-1)^{m-1-i}\binom{m}{i}t^i\in\mathbb{Z}[t].$$
Hence the degree of $c_m(t)$ is $m-1$ with the leading coefficient $m!$. For example, we have
\begin{align*}
&c_1(t)=1,\\
&c_2(t)=2t-1,\\
&c_3(t)=2(3t^2-3t+1),\\
&c_4(t)=6(4t^3-6t^2+4t-1).
\end{align*}
It is easy to see that
$$c_m(0)=(m-1)!(-1)^{m-1},\qquad c_m(1)=(m-1)!,$$
and if $m\geqslant 2$, it holds
\begin{align}
c_{m+1}(t)=m(2t-1)c_m(t)-m(m-1)(t^2-t)c_{m-1}(t).
\label{Eq:Sym-Coe-Recur}
\end{align}

For any $n\in\mathbb{N}$, denote by $\mathcal{P}_n$ the set of all partitions of the set $\{1,2,\ldots,n\}$. For a partition $\Pi=\{P_1,\ldots,P_i\}\in\mathcal{P}_n$ with $l_j=\# P_j$, define
$$c_{\Pi}(t)=\prod\limits_{j=1}^ic_{l_j}(t)\in\mathbb{Z}[t].$$
For example, if $\Pi=\{\{1,2\},\{3,4\},\{5\}\}\in\mathcal{P}_5$, we have
$$c_{\Pi}(t)=c_2(t)^2c_1(t)=(2t-1)^2.$$

Let $S_n$ be the symmetric group of degree $n$. Now we have the following identity in $\mathfrak{h}^1_t$, which is a generalization of \cite[Lemma 5.1]{Li-Qin-2} because for any $\Pi=\{P_1,\ldots,P_i\}\in\mathcal{P}_n$ with $l_j=\# P_j$, we have
$$c_{\Pi}(0)=(-1)^{n-i}\prod\limits_{j=1}^i(l_j-1)!,\quad c_{\Pi}(1)=\prod\limits_{j=1}^i(l_j-1)!.$$

\begin{lem}\label{Lem:SymSum}
Let $\mathbf{k}=(k_1,\ldots,k_n)$ be an index. We have
\begin{align}
\sum\limits_{\sigma\in S_n}z_{k_{\sigma(1)}}\cdots z_{k_{\sigma(n)}}=\sum\limits_{\Pi=\{P_1,\ldots,P_i\}\in\mathcal{P}_n}c_{\Pi}(t)z_{\mathbf{k},P_1}\tast\cdots\tast z_{\mathbf{k},P_i},
\label{Eq:SymSum}
\end{align}
where  $z_{\mathbf{k},P_j}=z_{\sum\limits_{l\in P_j}k_l}$.
\end{lem}

\proof
We proceed by induction on $n$, which is similar as the proof of \cite[Lemma 5.1]{Li-Qin-2}. The case of $n=1$ is trivial, and the case of $n=2$ follows from the fact
$$z_{k_1}\tast z_{k_2}=z_{k_1}z_{k_2}+z_{k_2}z_{k_1}+(1-2t)z_{k_1+k_2}.$$
Now assume that $n\geqslant 2$, $\mathbf{k}=(k_1,\ldots,k_n)$ and $\mathbf{k}'=(k_1,\ldots,k_n,k_{n+1})$. As
\begin{align*}
&z_{k_{n+1}}\tast \sum\limits_{\sigma\in S_n}z_{k_{\sigma(1)}}\cdots z_{k_{\sigma(n)}}
=\sum\limits_{\sigma\in S_n}\sum\limits_{j=1}^{n+1}z_{k_{\sigma(1)}}\cdots z_{k_{\sigma(j-1)}}z_{k_{n+1}}z_{k_{\sigma(j)}}\cdots z_{k_{\sigma(n)}}\\
&\qquad+(1-2t)\sum\limits_{\sigma\in S_n}\sum\limits_{j=1}^nz_{k_{\sigma(1)}}\cdots z_{k_{\sigma(j-1)}}z_{k_{\sigma(j)}+k_{n+1}}z_{k_{\sigma(j+1)}}\cdots z_{k_{\sigma(n)}}\\
&\qquad+(t^2-t)\sum\limits_{\sigma\in S_n}\sum\limits_{j=1}^{n-1}z_{k_{\sigma(1)}}\cdots z_{k_{\sigma(j-1)}}z_{k_{\sigma(j)}+k_{\sigma(j+1)}+k_{n+1}}z_{k_{\sigma(j+2)}}\cdots z_{k_{\sigma(n)}},
\end{align*}
we obtain
\begin{align*}
&\sum\limits_{\sigma\in S_{n+1}}z_{k_{\sigma(1)}}\cdots z_{k_{\sigma(n+1)}}=z_{k_{n+1}}\tast \sum\limits_{\sigma\in S_n}z_{k_{\sigma(1)}}\cdots z_{k_{\sigma(n)}}\\
&\quad+(2t-1)\sum\limits_{j=1}^n\sum\limits_{\sigma\in S_n}z_{k_{\sigma(1)}^{(j)}}\cdots z_{k_{\sigma(n)}^{(j)}}-2(t^2-t)\sum\limits_{1\leqslant i<j\leqslant n}\sum\limits_{\sigma\in S_{n-1}}z_{k_{\sigma(1)}^{(i,j)}}\cdots z_{k_{\sigma(n-1)}^{(i,j)}},
\end{align*}
where
$$\mathbf{k}^{(j)}=(k_1,\ldots,k_{j-1},k_j+k_{n+1},k_{j+1},\ldots,k_n)=(k_1^{(j)},\ldots,k_n^{(j)})$$
and
$$\mathbf{k}^{(i,j)}=(k_i+k_j+k_{n+1},k_1,\ldots,\breve{k_i},\ldots,\breve{k_j},\ldots,k_n)=(k_1^{(i,j)},\ldots,k_{n-1}^{(i,j)}).$$
By inductive hypothesis on $\mathbf{k}$, $\mathbf{k}^{(j)}$ and $\mathbf{k}^{(i,j)}$, we get
\begin{align*}
&\sum\limits_{\sigma\in S_{n+1}}z_{k_{\sigma(1)}}\cdots z_{k_{\sigma(n+1)}}\\
=&\sum\limits_{\Pi=\{P_1,\ldots,P_r\}\in \mathcal{P}_n}c_{\Pi}(t)z_{\mathbf{k},P_1}\tast\cdots\tast z_{\mathbf{k},P_r}\tast z_{k_{n+1}}\\
&+(2t-1)\sum\limits_{j=1}^n\sum\limits_{\Pi=\{P_1,\ldots,P_r\}\in \mathcal{P}_n}c_{\Pi}(t)z_{\mathbf{k}^{(j)},P_1}\tast\cdots\tast z_{\mathbf{k}^{(j)},P_r}\\
&-2(t^2-t)\sum\limits_{1\leqslant i<j\leqslant n}\sum\limits_{\Pi=\{P_1,\ldots,P_r\}\in \mathcal{P}_{n-1}}c_{\Pi}(t)z_{\mathbf{k}^{(i,j)},P_1}\tast\cdots\tast z_{\mathbf{k}^{(i,j)},P_r}\\
=&\sum\limits_{\Pi=\{P_1,\ldots,P_r\}\in \mathcal{P}_{n+1}}\tilde{c}_{\Pi}(t)z_{\mathbf{k}',P_1}\tast\cdots\tast z_{\mathbf{k}',P_r}.
\end{align*}
Here for $\Pi=\{P_1,\ldots,P_r\}\in\mathcal{P}_{n+1}$ with $P_r=\{n+1\}$, we have
$$\tilde{c}_{\Pi}(t)=c_{\Pi'}(t)=c_{\Pi}(t),$$
where $\Pi'=\{P_1,\ldots,P_{r-1}\}\in\mathcal{P}_n$. And for $\Pi=\{P_1,\ldots,P_r\}\in\mathcal{P}_{n+1}$ with $P_r=\{j<n+1\}$, we have
$$\tilde{c}_{\Pi}(t)=(2t-1)c_{\Pi'}(t)=c_{\Pi}(t),$$
where $\Pi'=\{P_1,\ldots,P_{r-1},\{j\}\}\in\mathcal{P}_n$. Finally, for $\Pi=\{P_1,\ldots,P_r\}\in\mathcal{P}_{n+1}$, with $n+1\in P_r$ and $\# P_r=m+1\geqslant 3$, set $l_j=\# P_j$ with $j=1,\ldots,r-1$. Then we get
$$\tilde{c}_{\Pi}(t)=(2t-1)mc_m(t)\prod\limits_{j=1}^{r-1}c_{l_{j}}(t)-2(t^2-t)\binom{m}{2}c_{m-1}(t)\prod\limits_{j=1}^{r-1}c_{l_{j}}(t).$$
Using \eqref{Eq:Sym-Coe-Recur}, we find
$$\tilde{c}_{\Pi}(t)=c_{m+1}(t)\prod\limits_{j=1}^{r-1}c_{l_{j}}(t)=c_{\Pi}(t).$$
Hence we have shown that for any $\Pi\in\mathcal{P}_{n+1}$, it holds $\tilde{c}_{\Pi}(t)=c_{\Pi}(t)$. At last, we get
$$\sum\limits_{\sigma\in S_{n+1}}z_{k_{\sigma(1)}}\cdots z_{k_{\sigma(n+1)}}=\sum\limits_{\Pi=\{P_1,\ldots,P_r\}\in\mathcal{P}_{n+1}}c_{\Pi}(t)z_{\mathbf{k}',P_1}\tast\cdots\tast z_{\mathbf{k}',P_i},$$
which completes the proof.
\qed

Applying the map $Z^t$ to \eqref{Eq:SymSum} under the condition that all $k_i\geqslant 2$, we obtain the following symmetric sum formula of $t$-MZVs.

\begin{thm}[Symmetric sum formula]
Let $\mathbf{k}=(k_1,\ldots,k_n)$ be an index with all $k_i\geqslant 2$. We have
\begin{align*}
\sum\limits_{\sigma\in S_n}\zeta^t(k_{\sigma(1)},\ldots, k_{\sigma(n)})=\sum\limits_{\Pi=\{P_1,\ldots,P_i\}\in\mathcal{P}_n}c_{\Pi}(t)\prod\limits_{j=1}^i\zeta\left(\sum\limits_{l\in P_j}k_l\right).
\end{align*}
\end{thm}


\subsection{Hoffman's relation}

Hoffman's relation of $t$-MZVs was introduced in \cite{Li-Qin2017,Wakabayashi}. Here we recall a proof of Hoffman's relation. In fact, for $w=z_{k_1}\cdots z_{k_n}\in\mathfrak{h}_t^1$ with $n,k_1,\ldots,k_n\in\mathbb{N}$, using \cite[Lemma 3.1]{Li-Qin}, we have
$$\partial_1^t(w)=S_{-t}\circ\partial_1\circ S_t(w)=S_{-t}\left[y\shuffle S_t(w)-y\ast S_t(w)\right]=y\tshuffle w-y\tast w.$$
Therefore from \cite{Li-Qin2017,Wakabayashi}, we get
\begin{align*}
\partial_1^t(w)=&y\tshuffle w-y\tast w=\sum\limits_{i=1}^n\sum\limits_{j=2}^{k_i}z_{k_1}\cdots z_{k_{i-1}}z_jz_{k_i+1-j}z_{k_{i+1}}\cdots z_{k_n}\\
&-\sum\limits_{i=1}^n\left[1+(k_i+\delta_{ni}-2)t\right]z_{k_1}\cdots z_{k_{i-1}}z_{k_i+1}z_{k_{i+1}}\cdots z_{k_n}\\
&+(t-t^2)\sum\limits_{i=1}^{n-1}z_{k_1}\cdots z_{k_{i-1}}z_{k_i+k_{i+1}+1}z_{k_{i+2}}\cdots z_{k_n},
\end{align*}
where $\delta_{ij}$ is Kronecker's delta symbol defined by
$$\delta_{ij}=\begin{cases}
1 & \text{\;if\;} i=j,\\
0 & \text{\;if\;} i\neq j.
\end{cases}$$
Hence if $k_1\geqslant 2$, we have
\begin{align*}
&\sum\limits_{i=1}^n\left[1+(k_i+\delta_{ni}-2)t\right]\zeta^t(k_1,\ldots,k_{i-1},k_i+1,k_{i+1},\ldots,k_n)\\
=&\sum\limits_{i=1}^n\sum\limits_{j=2}^{k_i}\zeta^t(k_1,\ldots,k_{i-1},j,k_i+1-j,k_{i+1},\ldots,k_n)\\
&+(t-t^2)\sum\limits_{i=1}^{n-1}\zeta^t(k_1,\ldots,k_{i-1},k_i+k_{i+1}+1,k_{i+2},\ldots,k_n),
\end{align*}
which was called Hoffman's relation in \cite{Li-Qin2017,Wakabayashi}.


\subsection{Sum formula}

The sum formula of $t$-MZVs was given by S. Yamamoto in \cite{Yamamoto} by using the sum formula of MZVs. Two different proofs of the sum formula were given in \cite{Li-Qin2017}, one proof is analytical and the other one is algebraic. We give an algebraic proof here. Let $k,n\in\mathbb{N}$ with the condition $k>n$. Then
$$\sigma_{n-1}^t(z_{k-n+1})=S_{-t}(\sigma_{n-1}(z_{k-n+1}))=S_{-t}(z_k)=z_k.$$
By \cite[Lemma 3.2]{Li-Qin}, we have
$$\overline{\sigma_{n-1}}(z_{k-n+1})=x(x^{k-n-1}\shuffle y^{n-1})y,$$
which induces that
$$\overline{\sigma_{n-1}}^t(z_{k-n+1})=S_{-t}(x(x^{k-n-1}\shuffle y^{n-1})y)=x\sigma_{-t}(x^{k-n-1}\shuffle y^{n-1})y.$$
As $\sigma_{-t}$ is an automorphism of the algebra $(\mathfrak{h}_t,\shuffle)$, we get
\begin{align*}
\overline{\sigma_{n-1}}^t(z_{k-n+1})=&x(x^{k-n-1}\shuffle (-tx+y)^{n-1})y\\
=&\sum\limits_{i=1}^nx(x^{k-n-1}\shuffle x^{n-i}\shuffle y^{i-1})y(-t)^{n-i}\\
=&\sum\limits_{i=1}^n\binom{k-i-1}{n-i}x(x^{k-i-1}\shuffle y^{i-1})y(-t)^{n-i}.
\end{align*}
Hence we have
$$\left(\overline{\sigma_{n-1}}^t-\sigma_{n-1}^t\right)(z_{k-n+1})=\sum\limits_{i=1}^n\binom{k-i-1}{n-i}x(x^{k-i-1}\shuffle y^{i-1})y(-t)^{n-i}-z_k.$$
Then as in \cite[Appendix A]{Li-Qin2017} or as in \cite{Li-Qin}, we get
$$\sum\limits_{\mathbf{k}\in I^0_{k,n}}\zeta^t(\mathbf{k})=\left[\sum\limits_{i=1}^n\binom{k-i-1}{n-i}t^{n-i}\right]\zeta(k),$$
which is the sum formula of $t$-MZVs and is equivalent to \eqref{Eq:SumFormula}.


\subsection{Shuffle regularized sum formula}

In \cite[Theorem 1.2]{Kaneko-Sataka}, M. Kaneko and M. Sakata proved that for any $k,n\in\mathbb{N}$, it holds
\begin{align}
Z\left(\reg_{\shuffle}\left(\sum\limits_{w\in A^{\ast}\cap\mathfrak{h}^1\atop d_x(w)=k,d_y(w)=n}w\right)\right)=(-1)^{n-1}\zeta(k+1,\{1\}^{n-1}),
\label{Eq:ShuRegSum-MZV}
\end{align}
which was called the shuffle regularized sum formula of MZVs there.  Here and below, let $\{m_1,\ldots,m_l\}^n$ be $n$ repetitions of $m_1,\ldots,m_l$. In fact, \eqref{Eq:ShuRegSum-MZV} is also an immediate consequence of \cite[Lemma 3.3]{Li2010}, which claims that
\begin{align}
\reg_{\shuffle}\left(\sum\limits_{w\in A^{\ast}\cap\mathfrak{h}^1\atop d_x(w)=k,d_y(w)=n}w\right)=(-1)^{n-1}x^ky^n.
\label{Eq:ShuRegSum}
\end{align}
The following lemma can be regarded as a $t$-shuffle generalization of \eqref{Eq:ShuRegSum}.

\begin{lem}\label{Lem:t-Shu-Reg}
For any $k,n\in\mathbb{N}$, we have
\begin{align}
&\reg_{\tshuffle}\left(\sum\limits_{w\in A^{\ast}\cap\mathfrak{h}^1\atop d_x(w)=k,d_y(w)=n}w\right)=\sum\limits_{i=1}^{n}(-1)^{i-1}t^{n-i}\sum\limits_{k_1+\cdots+k_{i}=k+n\atop k_1\geqslant k+1,k_2,\ldots,k_{i}\geqslant 1}\binom{k_1}{k+1}z_{k_1}\cdots z_{k_{i}}.
\label{Eq:tShuRegSum}
\end{align}
\end{lem}

\proof
The left-hand side of \eqref{Eq:tShuRegSum} is
$$\reg_{\tshuffle}((x^k\shuffle y^{n-1})y)=S_{-t}(\reg_{\shuffle}(\sigma_t(x^{k}\shuffle y^{n-1})y)).$$
As $\sigma_t$ is an automorphism of the algebra $(\mathfrak{h}_t,\shuffle)$, the left-hand side of \eqref{Eq:tShuRegSum} becomes
\begin{align*}
&S_{-t}(\reg_{\shuffle}(x^k\shuffle(tx+y)^{n-1})y)=\sum\limits_{j=1}^{n}S_{-t}(\reg_{\shuffle}(x^k\shuffle x^{n-j}\shuffle y^{j-1})y)t^{n-j}\\
=&\sum\limits_{j=1}^{n}\binom{k+n-j}{k}S_{-t}(\reg_{\shuffle}(x^{k+n-j}\shuffle y^{j-1})y)t^{n-j},
\end{align*}
which by \eqref{Eq:ShuRegSum} is
$$\sum\limits_{j=1}^{n}\binom{k+n-j}{k}(-1)^{j-1}S_{-t}(x^{k+n-j}y^{j})t^{n-j}.$$
Hence the left-hand side of \eqref{Eq:tShuRegSum} equals
\begin{align*}
&\sum\limits_{j=1}^{n}\binom{k+n-j}{k}(-1)^{j-1}x^{k+n-j}(-tx+y)^{j-1}yt^{n-j}\\
=&\sum\limits_{j=1}^{n}\sum\limits_{i=1}^{j}\binom{k+n-j}{k}(-1)^{i-1}x^{k+n-j}(x^{j-i}\shuffle y^{i-1})yt^{n-i}\\
=&\sum\limits_{i=1}^n\sum\limits_{j=i}^n\binom{k+n-j}{k}(-1)^{i-1}\sum\limits_{k_1+\cdots+k_i=k+n\atop k_1\geqslant k+n-j+1,k_2,\ldots,k_i\geqslant 1}z_{k_1}\cdots z_{k_i}t^{n-i}\\
=&\sum\limits_{i=1}^n(-1)^{i-1}t^{n-i}\sum\limits_{k_1+\cdots+k_i=k+n\atop k_1\geqslant k+1,k_2,\ldots,k_i\geqslant 1}z_{k_1}\cdots z_{k_i}\sum\limits_{j=k+n+1-k_1}^n\binom{k+n-j}{k},
\end{align*}
which is just the right-hand side of \eqref{Eq:tShuRegSum}.
\qed

From Lemma \ref{Lem:t-Shu-Reg}, we get the shuffle regularized sum formula of $t$-MZVs, which generalizes \eqref{Eq:ShuRegSum-MZV}.

\begin{thm}[Shuffle regularized sum formula]
For any $k,n\in\mathbb{N}$, we have
\begin{align*}
&Z^t\left(\reg_{\tshuffle}\left(\sum\limits_{w\in A^{\ast}\cap\mathfrak{h}^1\atop d_x(w)=k,d_y(w)=n}w\right)\right)\\
=&\sum\limits_{i=1}^{n}(-1)^{i-1}t^{n-i}\sum\limits_{k_1+\cdots+k_{i}=k+n\atop k_1\geqslant k+1,k_2,\ldots,k_{i}\geqslant 1}\binom{k_1}{k+1}\zeta^t(k_1,\ldots,k_{i}).
\end{align*}
\end{thm}


\subsection{A formula for height one $t$-MZVs}

For an index $\mathbf{k}=(k_1,\ldots,k_n)$, we define its height by
$$\height(\mathbf{k})=\#\{j\mid 1\leqslant j\leqslant n,k_j\geqslant 2\},$$
which was first introduced by Y. Ohno and D. Zagier in \cite{Ohno-Zagier}. For example, an admissible index $\mathbf{k}$ is of height one if and only if $\mathbf{k}$ has the form
$$\mathbf{k}=(k+1,\{1\}^{l-1}),\quad (\exists k,l\in\mathbb{N}).$$

M. Kaneko and M. Sakata gave a formula for height one MZVs in \cite[Theorem 1.1]{Kaneko-Sataka}. They showed that for any $k,l\in\mathbb{N}$, it holds
\begin{align}
\zeta(k+1,\{1\}^{l-1})=\sum\limits_{j=1}^{\min\{k,l\}}(-1)^{j-1}\sum\limits_{\wt(\mathbf{k})=k,\wt(\mathbf{l})=l\atop \dep(\mathbf{k})=\dep(\mathbf{l})=j}\zeta(\mathbf{k}+\mathbf{l}),
\label{Eq:HeightOne-MZV}
\end{align}
where for $\mathbf{k}=(k_1,\ldots,k_j)$ and $\mathbf{l}=(l_1,\ldots,l_j)$, $\mathbf{k}+\mathbf{l}=(k_1+l_1,\ldots,k_j+l_j)$. We remark that the above formula is also a consequence of the identity
\begin{align}
\Delta_v\left(\frac{1}{1-xu-xv+x(x+y)uv}xy\right)=\frac{x}{1-xu}\frac{y}{1-yv},
\label{Eq:HeightOne-Func}
\end{align}
where $u$ and $v$ are variables commuting with each other as well as $x$ and $y$,  and $\Delta_v$ is defined in \cite{Ihara-Kaneko-Zagier} by
$$\Delta_v=\exp\left(\sum\limits_{n=1}^\infty\frac{\partial_n}{n}v^n\right).$$
Here \eqref{Eq:HeightOne-Func} is essential the last equation of \cite{Kajikawa}. And to get \eqref{Eq:HeightOne-MZV} from \eqref{Eq:HeightOne-Func}, one can use the computation
\begin{align*}
&\frac{1}{1-xu-xv+x(x+y)uv}xy=\frac{1}{1+\frac{xu}{1-xu}\frac{1}{1-xv}yv}\frac{xu}{1-xu}\frac{1}{1-xv}yvu^{-1}v^{-1}\\
=&\sum\limits_{j=1}^\infty(-1)^{j-1}\left(\frac{xu}{1-xu}\frac{1}{1-xv}yv\right)^ju^{-1}v^{-1}\\
=&\sum\limits_{j=1}^{\infty}(-1)^{j-1}\sum\limits_{k_1,\ldots,k_j\geqslant 1\atop l_1,\ldots,l_j\geqslant 1}x^{k_1+l_1-1}y\cdots x^{k_j+l_j-1}yu^{k_1+\cdots+k_j-1}v^{l_1+\cdots+l_j-1}\\
=&\sum\limits_{k,l\geqslant 1}\left(\sum\limits_{j=1}^{\min\{k,l\}}(-1)^{j-1}\sum\limits_{{k_1+\cdots+k_j=k\atop l_1+\cdots+l_j=l}\atop k_1,\ldots,k_j,l_1,\ldots,l_j\geqslant 1}z_{k_1+l_1}\cdots z_{k_j+l_j}\right)u^{k-1}v^{l-1},
\end{align*}
and the derivation relations of MZVs introduced by K. Ihara, M. Kaneko and D. Zagier in \cite{Ihara-Kaneko-Zagier}.

We would like to generalize \eqref{Eq:HeightOne-Func} to $t$-MZVs. Define
$$\Delta_v^t=\exp\left(\sum\limits_{n=1}^\infty\frac{\partial_n^t}{n}v^n\right).$$
Then since $\partial_n^t=S_{-t}\circ \partial_n\circ S_t$, we get $\Delta_v^t=S_{-t}\circ \Delta_v\circ S_t$. Note that for any $w\in\mathfrak{h}_t^0$, we have
\begin{align}
Z^t(\Delta_v^t(w))=Z^t(w),
\label{Eq:DerRel-tMZV}
\end{align}
which is from the extended double shuffle relations of $t$-MZVs.

Now a generalization of \eqref{Eq:HeightOne-Func} is displayed in the following lemma.

\begin{lem}\label{Lem:t-HeightOne-Func}
In $\mathfrak{h}_t^0[[u,v]]$, we have
\begin{align}
\Delta_v^t\left(\frac{1}{1-xu-xv+x[(1-t)x+y]uv}\frac{x}{1-txv}y\right)=\frac{x}{1-xu}\frac{y}{1-yv}.
\label{Eq:t-HeightOne-Func}
\end{align}
\end{lem}

\proof
The left-hand side of \eqref{Eq:t-HeightOne-Func} is
$$(S_{-t}\circ \Delta_v)\left(\frac{1}{1-xu-xv+x(x+y)uv}\frac{x}{1-txv}y\right).$$
Since
$$\Delta_v(x)=x\frac{1}{1-yv},\quad \Delta_v(y)=(1-xv-yv)\frac{y}{1-yv},$$
we get
$$\Delta_v\left(1-xu-xv+x(x+y)uv\right)=x\frac{1}{1-yv}(1-xv-yv)\frac{1-xu}{x},$$
which implies that the left-hand side of \eqref{Eq:t-HeightOne-Func} is
$$S_{-t}\left(\frac{x}{1-xu}\frac{1}{1-xv-yv}(1-yv)\frac{1}{1-txv-yv}(1-xv-yv)\frac{y}{1-yv}\right).$$
Direct calculation shows that
$$(1-yv)\frac{1}{1-txv-yv}(1-xv-yv)=(1-xv-yv)\frac{1}{1-txv-yv}(1-yv),$$
which implies that the left-hand side of \eqref{Eq:t-HeightOne-Func} becomes
$$S_{-t}\left(\frac{x}{1-xu}\frac{1}{1-txv-yv}y\right).$$
Now it is easy to finish the proof.
\qed

From Lemma \ref{Lem:t-HeightOne-Func}, we get some formulas for the height one $t$-MZVs which generalize \eqref{Eq:HeightOne-MZV}.

\begin{thm}
For any $k,l\in\mathbb{N}$, we have
\begin{align}
\zeta^t(k+1,\{1\}^{l-1})=&\sum\limits_{j=1}^{\min\{k,l\}}(-1)^{j-1}\sum\limits_{{k_1+\cdots+k_j=k\atop l_1+\cdots+l_j=l}\atop k_1,\ldots,k_j,l_1,\ldots,l_j\geqslant 1}\frac{1-t^{l_j}}{1-t}
\nonumber\\
&\times Z^t(x^{k_1+l_1-1}(-tx+y)\cdots x^{k_{j-1}+l_{j-1}-1}(-tx+y)x^{k_j+l_j-1}y),
\label{Eq:HeightOne-tMZV}
\end{align}
or equivalently
\begin{align}
\zeta^t(k+1,\{1\}^{l-1})=&\sum\limits_{j=1}^{\min\{k,l\}}(-1)^{j-1}\sum\limits_{{k_1+\cdots+k_j=k\atop l_1+\cdots+l_j=l}\atop k_1,\ldots,k_j,l_1,\ldots,l_j\geqslant 1}\frac{1-t^{l_j}}{1-t}\zeta(k_1+l_1,\ldots,k_j+l_j).
\label{Eq:HeightOne-tMZV-1}
\end{align}
In particular, we have
\begin{align}
\zeta^{\star}(k+1,\{1\}^{l-1})=\sum\limits_{j=1}^{\min\{k,l\}}(-1)^{j-1}\sum\limits_{{k_1+\cdots+k_j=k\atop l_1+\cdots+l_j=l}\atop k_1,\ldots,k_j,l_1,\ldots,l_j\geqslant 1}l_j\zeta(k_1+l_1,\ldots,k_j+l_j).
\label{Eq:HeightOne-MZSV}
\end{align}
\end{thm}

\proof
Since
\begin{align*}
&\frac{1}{1-xu-xv+x[(1-t)x+y]uv}\frac{x}{1-txv}y\\
=&\frac{1}{1+\frac{xu}{1-xu}\frac{1}{1-xv}(-tx+y)v}\frac{xu}{1-xu}\frac{1}{1-xv}\frac{1}{1-txv}yu^{-1}\\
=&\sum\limits_{j=1}^\infty(-1)^{j-1}\sum\limits_{{k_1,\ldots,k_j\geqslant 1\atop l_1,\ldots,l_j\geqslant 1}\atop i\geqslant 0}x^{k_1+l_1-1}(-tx+y)\cdots x^{k_{j-1}+l_{j-1}-1}(-tx+y)\\
&\quad \times x^{k_j+l_j+i-1}yt^iu^{k_1+\cdots+k_i-1}v^{l_1+\cdots+l_j+i-1},
\end{align*}
then by \eqref{Eq:DerRel-tMZV}, we have
\begin{align*}
&\zeta^t(k+1,\{1\}^{l-1})=\sum\limits_{j=1}^{\min\{k,l\}}(-1)^{j-1}\sum\limits_{{k_1+\cdots+k_j=k\atop l_1+\cdots+l_j+i=l}\atop k_1,\ldots,k_j,l_1,\ldots,l_j\geqslant 1,i\geqslant 0}\\
&\qquad\times Z^t(x^{k_1+l_1-1}(-tx+y)\cdots x^{k_{j-1}+l_{j-1}-1}(-tx+y)x^{k_j+l_j+i-1}y)t^i,
\end{align*}
which implies \eqref{Eq:HeightOne-tMZV}. We get \eqref{Eq:HeightOne-tMZV-1} from \eqref{Eq:HeightOne-tMZV} by using the fact $Z^t=Z\circ S_t$, and get \eqref{Eq:HeightOne-MZSV} from \eqref{Eq:HeightOne-tMZV} by setting $t=1$.
\qed

As in \cite{Kaneko-Sataka}, \eqref{Eq:HeightOne-tMZV} is nothing but the sum formula of double $t$-MZVs in the case of $l=2$.


\subsection{A weighted sum formula}

In \cite{Guo-Xie2009}, L. Guo and B. Xie provided a weighted sum formula of MZVs. And in \cite{Li-Qin}, C. Qin and the author simplified the proof of this weighted sum formula and gave its MZSVs version. In this subsection, we shall generalize the weighted sum formulas to $t$-MZVs. The method used here is similar as that in \cite{Guo-Xie2009,Li-Qin}.

To save spaces, we use the following notation. If $\mathbf{k}\in I_{k,n}$, we always mean that $\mathbf{k}=(k_1,\ldots,k_n)$ unless stated otherwise, and set $z_{\mathbf{k}}=z_{k_1}\cdots z_{k_n}$.

\begin{lem}
For any integers $k,n$ with $k>n\geqslant 2$, we have
\begin{align}
&\sum\limits_{l+k_1+\cdots+k_{n-1}=k\atop l,k_i\geqslant 1,k_1\geqslant 2}z_l\tast z_{k_1}\cdots z_{k_{n-1}}=n\sum\limits_{\mathbf{k}\in I_{k,n}^0}z_{\mathbf{k}}+\sum\limits_{\mathbf{k}\in I_{k,n},k_1=1}z_{\mathbf{k}}-\sum\limits_{\mathbf{k}\in I_{k,n},k_2=1}z_{\mathbf{k}}\nonumber\\
&\qquad\qquad\qquad+(1-2t)(k-n)\sum\limits_{\mathbf{k}\in I^0_{k,n-1}}z_{\mathbf{k}}+(t^2-t)\sum\limits_{\mathbf{k}\in I^0_{k,n-2}}b_{\mathbf{k}}z_{\mathbf{k}},
\label{Eq:Sum-tast}
\end{align}
where $I^0_{k,n-2}=\varnothing$ for $n=2$, and for $n>2$ and $\mathbf{k}\in I_{k,n-2}^0$,
$$b_{\mathbf{k}}=\binom{k_1-2}{2}+\binom{k_2-1}{2}+\cdots+\binom{k_{n-2}-1}{2}.$$
\end{lem}

\proof
The proof is similar as that in \cite[Lemmas 3.3, 3.4]{Li-Qin}. By the definition of the product $\tast$, we have
\begin{align}
&z_l\tast z_{k_1}\cdots z_{k_{n-1}}=\sum\limits_{i=0}^{n-1}z_{k_1}\cdots z_{k_i}z_lz_{k_{i+1}}\cdots z_{k_{n-1}}\nonumber\\
&\qquad\qquad+(1-2t)\sum\limits_{i=1}^{n-1}z_{k_1}\cdots z_{k_{i-1}}z_{l+k_i}z_{k_{i+1}}\cdots z_{k_{n-1}}\nonumber\\
&\qquad\qquad+(t^2-t)\sum\limits_{i=1}^{n-2}z_{k_1}\cdots z_{k_{i-1}}z_{l+k_i+k_{i+1}}z_{k_{i+2}}\cdots z_{k_{n-1}}.
\label{Eq:tast-1-n}
\end{align}
Therefore we find the left-hand side of \eqref{Eq:Sum-tast} is
\begin{align*}
&\sum\limits_{\mathbf{k}\in I_{k,n},k_2\geqslant 2}z_{\mathbf{k}}+(n-1)\sum\limits_{\mathbf{k}\in I^0_{k,n}}z_{\mathbf{k}}+(1-2t)\sum\limits_{\mathbf{k}\in I^0_{k,n-1}}(k_1-2)z_{\mathbf{k}}\\
&+(1-2t)\sum\limits_{i=2}^{n-1}\sum\limits_{\mathbf{k}\in I^0_{k,n-1}}(k_i-1)z_{\mathbf{k}}+(t^2-t)\sum\limits_{\mathbf{k}\in I^0_{k,n-2}}\binom{k_1-2}{2}z_{\mathbf{k}}\\
&+(t^2-t)\sum\limits_{i=2}^{n-2}\sum\limits_{\mathbf{k}\in I^0_{k,n-2}}\binom{k_i-1}{2}z_{\mathbf{k}}.
\end{align*}
Now it is easy to finish the proof.
\qed

To deal with the product $\tshuffle$, we first need a formula similar as \eqref{Eq:tast-1-n}. The formula is a generalization of that for shuffle product of MZVs appearing in the proof of \cite[Lemma 3.3]{Li-Qin} and of MZSVs appearing in the proof of \cite[Lemma 3.4]{Li-Qin}. See also \cite[Proposition 2.7]{Li-Qin-shuffle}.

\begin{lem}
For any $n,l,k_1,\ldots,k_{n-1}\in\mathbb{N}$ with $n\geqslant 2$, we have
\begin{align}
&z_l\tshuffle z_{k_1}\cdots z_{k_{n-1}}=\sum\limits_{i=1}^{n-1}\sum\limits_{\alpha_1+\cdots+\alpha_{i+1}\atop =l+k_1+\cdots+k_{i},\alpha_j\geqslant 1}\prod\limits_{j=1}^{i-1}\binom{\alpha_j-1}{k_j-1}\binom{\alpha_i-1}{k_i-\alpha_{i+1}}\nonumber\\
&\qquad\times (z_{\alpha_1}\cdots z_{\alpha_{i+1}}z_{k_{i+1}}\cdots z_{k_{n-1}}-tz_{\alpha_1}\cdots z_{\alpha_{i-1}}z_{\alpha_i+\alpha_{i+1}}z_{k_{i+1}}\cdots z_{k_{n-1}})\nonumber\\
&\qquad+\sum\limits_{\alpha_1+\cdots+\alpha_n\atop =l+k_1+\cdots+k_{n-1},\alpha_j\geqslant 1}\prod\limits_{j=1}^{n-1}\binom{\alpha_j-1}{k_j-1}(z_{\alpha_1}\cdots z_{\alpha_n}-tz_{\alpha_1}\cdots z_{\alpha_{n-2}}z_{\alpha_{n-1}+\alpha_n}).
\label{Eq:tshuffle-1-n}
\end{align}
\end{lem}

\proof
We use the combinatorial description of the shuffle product to prove this formula. For more applications of the combinatorial description of the shuffle product, one may refer to \cite{Li-Qin-shuffle}. To compute the $t$-shuffle product, we consider the $y$'s first. Denote the letter $y$ in $z_l$ by $Y$. Then the shuffle product $z_l\shuffle z_{k_1}\cdots z_{k_{n-1}}$ is a linear combination of the following words:
\begin{itemize}
  \item [(i)] $x^{\alpha_1-1}y\cdots x^{\alpha_{i-1}-1}yx^{\alpha_i-1}Yx^{\alpha_{i+1}-1}y\cdots x^{\alpha_n-1}y$,\quad  ($1\leqslant i\leqslant n-1$),
  \item [(ii)] $x^{\alpha_1-1}y\cdots x^{\alpha_{n-1}-1}yx^{\alpha_n-1}Y$,
\end{itemize}
where $\alpha_1,\ldots,\alpha_n\in\mathbb{N}$ and satisfy the condition $\alpha_1+\cdots+\alpha_n=l+k_1+\cdots+k_{n-1}$.

Now we consider the $x$'s. For a word in the case of (i), the sources of $x$'s in $x^{\alpha_j-1}$ are
\begin{itemize}
  \item [(a)] $k_j-1$ times from $z_{k_j}$ and others from $z_l$ for $1\leqslant j\leqslant i-1$,
  \item [(b)] $k_i-\alpha_{i+1}$ times from $z_{k_i}$ and others from $z_l$ for $j=i$,
  \item [(c)] all from $z_{k_{j-1}}$ for $i+1\leqslant j\leqslant n$.
\end{itemize}
Hence we have $\alpha_{j}=k_{j-1}$ for $j=i+2,\ldots,n$, and the coefficient of this word is
$$\prod\limits_{j=1}^{i-1}\binom{\alpha_j-1}{k_j-1}\binom{\alpha_i-1}{k_i-\alpha_{i+1}}.$$
Different from the shuffle product, the $t$-shuffle product will have a new word coming from a word in the case of (i) by replacing $Y$ by $-tx$ with the same coefficient. Hence we get the first term of the right-hand side of \eqref{Eq:tshuffle-1-n}.

Similarly, one get the second term of the right-hand side of \eqref{Eq:tshuffle-1-n} after considering the words  in the case of (ii).
\qed

As in \cite{Guo-Xie2009}, for an index $\mathbf{k}=(k_1,k_2,\ldots,k_n)$, set
$$\mathcal{C}(\mathbf{k})=\mathcal{C}(k_1,\ldots,k_{n})=\sum\limits_{j=1}^{n}2^{k_1+\cdots+k_{j}-j}+2^{k_1+\cdots+k_{n}-n}\in\mathbb{N}.$$
We also set
$$\widetilde{\mathcal{C}}(\mathbf{k})=\mathcal{C}(k_1,\ldots,k_{n})-\mathcal{C}(k_2,\ldots,k_n),$$
where if $n=1$, $\mathcal{C}(k_2,\ldots,k_n)$ is treated as $1$.

From \eqref{Eq:tshuffle-1-n}, we get the $t$-shuffle version of \eqref{Eq:Sum-tast}.

\begin{lem}
For any integers $k,n$ with $k>n\geqslant 2$, we have
\begin{align}
&\sum\limits_{l+k_1+\cdots+k_{n-1}=k\atop l,k_i\geqslant 1,k_1\geqslant 2}z_l\tshuffle z_{k_1}\cdots z_{k_{n-1}}=-\sum\limits_{\mathbf{k}\in I_{k,n},k_2=1}z_{\mathbf{k}}+t\sum\limits_{\mathbf{k}\in I^0_{k,n-1}}z_{\mathbf{k}}\nonumber\\
&\qquad\qquad\qquad+\sum\limits_{\mathbf{k}\in I_{k,n}}\widetilde{\mathcal{C}}(\widetilde{\mathbf{k}})z_{\mathbf{k}}-t\sum\limits_{\mathbf{k}\in I_{k,n-1}}\left(\widetilde{\mathcal{C}}(\mathbf{k})-\widetilde{\mathcal{C}}(\widetilde{\mathbf{k}})\right)z_{\mathbf{k}},
\label{Eq:Sum-tshuffle}
\end{align}
where $\widetilde{\mathbf{k}}=(k_1,\ldots,k_{l-1})$ for any index $\mathbf{k}=(k_1,\ldots,k_l)$ with $l\geqslant 2$, and $\widetilde{\mathcal{C}}(\widetilde{\mathbf{k}})=k$ for $\mathbf{k}=(k)$.
\end{lem}

\proof
The proof is similar as that in the special case of $t=1$, which can be found in the proof of \cite[Lemma 3.4]{Li-Qin}.
\qed

Using \eqref{Eq:Sum-tast} and \eqref{Eq:Sum-tshuffle}, we immediately get the following result.

\begin{cor}
For any integers $k,n$ with $k>n\geqslant 2$, we have
\begin{align}
&\sum\limits_{l+k_1+\cdots+k_{n-1}=k\atop l,k_i\geqslant 1,k_1\geqslant 2}\left(z_l\tshuffle z_{k_1}\cdots z_{k_{n-1}}-z_l\tast z_{k_1}\cdots z_{k_{n-1}}\right)\nonumber\\
=&\sum\limits_{\mathbf{k}\in I^0_{k,n}}\widetilde{\mathcal{C}}(\widetilde{\mathbf{k}})z_{\mathbf{k}}-n\sum\limits_{\mathbf{k}\in I^0_{k,n}}z_{\mathbf{k}}-t\sum\limits_{\mathbf{k}\in I^0_{k,n-1}}\left(\widetilde{\mathcal{C}}(\mathbf{k})-\widetilde{\mathcal{C}}(\widetilde{\mathbf{k}})\right)z_{\mathbf{k}}\nonumber\\
&+[t+(2t-1)(k-n)]\sum\limits_{\mathbf{k}\in I^0_{k,n-1}}z_{\mathbf{k}}+(t-t^2)\sum\limits_{\mathbf{k}\in I^0_{k,n-2}}b_{\mathbf{k}}z_{\mathbf{k}}.
\label{Eq:Sum-tshuffle-tast}
\end{align}
\end{cor}

Using the extended double shuffle relations and the sum formula \eqref{Eq:SumFormula}, we finally get a weighted sum formula of $t$-MZVs.

\begin{thm}[Weighted sum formula]
For any integers $k,n$ with $k>n\geqslant 2$, we have
\begin{align}
&\sum\limits_{\mathbf{k}\in I^0_{k,n}}\widetilde{\mathcal{C}}(\widetilde{\mathbf{k}})\zeta^t(\mathbf{k})-t\sum\limits_{\mathbf{k}\in I^0_{k,n-1}}\left(\widetilde{\mathcal{C}}(\mathbf{k})-\widetilde{\mathcal{C}}(\widetilde{\mathbf{k}})\right)\zeta^t(\mathbf{k})+(t-t^2)\sum\limits_{\mathbf{k}\in I^0_{k,n-2}}b_{\mathbf{k}}\zeta^t(\mathbf{k})\nonumber\\
=&\left\{\sum\limits_{i=1}^{n-2}\left[k\binom{k-1}{i}-(k-n+1)\binom{k-1}{i-1}\right]t^i(1-t)^{n-1-i}\right.\nonumber\\
&\qquad\left.+k(1-t)^{n-1}+\binom{k-1}{n-1}t^{n-1}\right\}\zeta(k),
\label{Eq:Weighted-Sum}
\end{align}
where $\widetilde{\mathbf{k}}=(k_1,\ldots,k_{l-1})$ for any index $\mathbf{k}=(k_1,\ldots,k_l)$ with $l\geqslant 2$, $\widetilde{\mathcal{C}}(\widetilde{\mathbf{k}})=k$ for $\mathbf{k}=(k)$, $I^0_{k,n-2}=\varnothing$ for $n=2$, and
$$b_{\mathbf{k}}=\binom{k_1-2}{2}+\binom{k_2-1}{2}+\cdots+\binom{k_{n-2}-1}{2}$$
for $n>2$ and $\mathbf{k}\in I_{k,n-2}^0$.
\end{thm}

Taking $t=0$ in \eqref{Eq:Weighted-Sum}, we get the weighted sum formula of MZVs which was first proved by L. Guo and B. Xie in \cite{Guo-Xie2009}. And taking $t=1$ in \eqref{Eq:Weighted-Sum}, we get the weighted sum formula of MZSVs in \cite{Li-Qin}. Setting $n=2$ in \eqref{Eq:Weighted-Sum} and using the sum formula for double $t$-MZVs, we get
$$\sum\limits_{k_1+k_2=k\atop k_1\geqslant 2,k_2\geqslant 1}2^{k_1}\zeta^t(k_1,k_2)=\left[k+1+(2^k-4)t\right]\zeta(k),\quad (k\geqslant 3).$$


\section{Some evaluation formulas}\label{Sec:EvaFormula}

In this section, we give some evaluation formulas of $t$-MZVs at even arguments, which are deduced from the extended double shuffle relations. As above, let $u$ and $v$ be variables commuting with each other as well as $x$ and $y$.


\subsection{Evaluation formulas of $\zeta^t(\{2k\}^n)$}

For any $m\in\mathbb{N}$, set
$$d_m(t)=\frac{1}{(m-1)!}c_m(t)=t^m-(t-1)^m=\sum\limits_{i=0}^{m-1}(-1)^{m-1-i}\binom{m}{i}t^i\in\mathbb{Z}[t].$$
For example, we have $d_1(t)=1$ and $d_2(t)=2t-1$. By \eqref{Eq:Sym-Coe-Recur}, for $m\geqslant 2$, it holds
\begin{align}
d_{m+1}(t)=(2t-1)d_m(t)-(t^2-t)d_{m-1}(t).
\label{Eq:tast-Recur}
\end{align}
Let $d_0(t)=0$, then \eqref{Eq:tast-Recur} is also valid for $m=1$.

The following lemma generalizes \cite[Lemma 3.10]{Li-Qin}.

\begin{lem}
For any $a,b,n\in\mathbb{N}$, we have
\begin{align}
\sum\limits_{m=1}^nd_m(t)z_{b+(m-1)a}\tast z_a^{n-m}=\sum\limits_{m=0}^{n-1}z_a^mz_bz_a^{n-1-m}.
\label{Eq:tast-Sum}
\end{align}
\end{lem}

\proof
We prove \eqref{Eq:tast-Sum} by induction on $n$. The case of $n=1$ is trivial. If $n=2$, the left-hand side of \eqref{Eq:tast-Sum} is
$$z_b\tast z_a+(2t-1)z_{b+a}=z_bz_a+z_az_b.$$
Thus \eqref{Eq:tast-Sum} holds for $n=2$.

Now assume that $n\geqslant 3$. Then the left-hand side of \eqref{Eq:tast-Sum} is
\begin{align*}
&d_n(t)z_{b+(n-1)a}+d_{n-1}(t)z_{b+(n-2)a}\tast z_a+\sum\limits_{m=1}^{n-2}d_m(t)\left[z_{b+(m-1)a}z_a^{n-m}\right.\\
&\left.+z_a(z_{b+(m-1)a}\tast z_a^{n-1-m})+(1-2t)z_{b+ma}z_a^{n-1-m}+(t^2-t)z_{b+(m+1)a}z_a^{n-2-m}\right].
\end{align*}
Using the inductive assumption, the above turns to
\begin{align*}
&\sum\limits_{m=1}^nd_m(t)z_{b+(m-1)a}z_a^{n-m}+\sum\limits_{m=1}^{n-1}(1-2t)d_m(t)z_{b+ma}z_a^{n-1-m}\\
&+\sum\limits_{m=1}^{n-2}(t^2-t)d_m(t)z_{b+(m+1)a}z_a^{n-2-m}+\sum\limits_{m=0}^{n-2}z_a^{m+1}z_bz_a^{n-2-m}.
\end{align*}
Finally, we finish the proof by \eqref{Eq:tast-Recur}.
\qed

Setting $a=b=k\in\mathbb{N}$ in \eqref{Eq:tast-Sum}, we get
\begin{align}
\sum\limits_{m=1}^nd_m(t)z_{mk}\tast z_k^{n-m}=nz_k^n, \quad (\forall n\in\mathbb{N}).
\label{Eq:tast-Sum-zk}
\end{align}
Then as a $t$-version of \cite[Corollary 2]{Ihara-Kaneko-Zagier}, we get the following result.

\begin{cor}
For any $k\in\mathbb{N}$, we have
\begin{align}
\frac{1}{1-z_ku}=\exp_{\tast}\left(\sum\limits_{n=1}^\infty\frac{1}{n}d_n(t)z_{nk}u^n\right).
\label{Eq:zkPower-exp}
\end{align}
In particular, if $k\geqslant 2$, we have
\begin{align}
1+\sum\limits_{n=1}^\infty\zeta^t(\{k\}^n)u^n=\exp\left(\sum\limits_{n=1}^\infty\frac{1}{n}d_n(t)\zeta(nk)u^n\right).
\label{Eq:zetat-zk-Power}
\end{align}
\end{cor}

\proof
Multiplying $u^n$ to both-sides of \eqref{Eq:tast-Sum-zk}, and then summing for $n$ from $1$ to $\infty$, we get
$$F(u)\tast\left(\sum\limits_{n=1}^\infty d_n(t)z_{nk}u^{n-1}\right)=F'(u),$$
where
$$F(u)=\frac{1}{1-z_ku}=\sum\limits_{n=0}^\infty z_k^nu^n.$$
Since $F(0)=1$, we obtain \eqref{Eq:zkPower-exp}.
\qed

In particular, since $d_n(0)=(-1)^{n-1}$, from \eqref{Eq:zetat-zk-Power} we get \cite[Corollary]{Ihara-Kaneko-Zagier}
$$1+\sum\limits_{n=1}^\infty\zeta(\{k\}^n)u^n=\exp\left(\sum\limits_{n=1}^\infty\frac{1}{n}(-1)^{n-1}\zeta(nk)u^n\right).$$
And since $d_n(1)=1$, we obtain
$$1+\sum\limits_{n=1}^\infty\zeta^{\star}(\{k\}^n)u^n=\exp\left(\sum\limits_{n=1}^\infty\frac{1}{n}\zeta(nk)u^n\right).$$

More generally, let $\mathfrak{z}$ be the $\mathbb{Q}[t]$-submodule of $\mathfrak{h}_t^1$ generated by $\{z_k\mid k\geqslant 1\}$. Define a circle product on $\mathfrak{z}$ by $\mathbb{Q}[t]$-bilinearity and the rule
$$z_k\circ z_l=z_{k+l}, \quad (\forall k,l\geqslant 1).$$
Then as \eqref{Eq:tast-Sum}, we find that for any $z\in \mathfrak{z}$ and any $n\in\mathbb{N}$, it holds
$$\sum\limits_{m=1}^nd_m(t)z^{\circ m}\tast z^{n-m}=nz^n,$$
which is a generalization of \eqref{Eq:tast-Sum-zk}. Hence as \eqref{Eq:zkPower-exp}, one gets
\begin{align*}
\frac{1}{1-zu}=&\exp_{\tast}\left(\sum\limits_{n=1}^\infty\frac{1}{n}d_n(t)z^{\circ n}u^n\right)\\
=&\exp_{\tast}\left\{\log_{\circ}\left(1-(t-1)zu\right)-\log_{\circ}(1-tzu)\right\},
\end{align*}
which generalizes \cite[Corollary 1]{Ihara-Kaneko-Zagier} and \cite[Proposition 3]{Ihara-Kajikawa-Ohno-Okuda}.

Using \eqref{Eq:zetat-zk-Power}, we can evaluate $\zeta^t(\{2k\}^n)$. For that purpose, we introduce a lemma. Let $\{B_n\}_{n=0}^\infty$ be the Bernoulli numbers defined by
$$\frac{u}{e^u-1}=\sum\limits_{n=0}^\infty\frac{B_n}{n!}u^n.$$

\begin{lem}
Let $k\in\mathbb{N}$. We have
\begin{align}
\sum\limits_{n=1}^\infty \frac{B_{2nk}}{2n(2nk)!}u^{2nk}=\sum\limits_{j=0}^{k-1}\left(\log\frac{e^{u\rho_k^j}-1}{u\rho_k^j}-\frac{1}{2}u\rho_k^j\right),
\label{Eq:Sum-BernoulliN}
\end{align}
where $\rho_k=e^{\frac{\pi\sqrt{-1}}{k}}$.
\end{lem}

\proof
Recall the formula
$$\log\frac{e^u-1}{u}-\frac{1}{2}u=\sum\limits_{n=1}^\infty \frac{B_{2n}}{2n(2n)!}u^{2n},$$
which was proved in \cite[Lemma 4.3]{Li2010}. Hence the right-hand side of \eqref{Eq:Sum-BernoulliN} becomes
$$\sum\limits_{n=1}^\infty\frac{B_{2n}}{2n(2n)!}u^{2n}\sum\limits_{j=0}^{k-1}\rho_k^{2nj}.$$
Since
$$\sum\limits_{j=0}^{k-1}\rho_k^{2nj}=\begin{cases}
0 & \text{if\;} k\nmid n,\\
k & \text{if\;} k\mid n,
\end{cases}$$
we get the result.
\qed

Furthermore, we need the evaluation formulas (see \cite[Eq. (3.26) and Eq. (3.27)]{Li-Qin})
\begin{align}
\zeta(\{2k\}^n)=&(-1)^n\sum\limits_{n_0+\cdots+n_{k-1}=nk\atop n_j\geqslant 0}\frac{\rho_k^{\sum\limits_{j=0}^{k-1}2jn_j}}{(2n_0+1)!\cdots(2n_{k-1}+1)!}\frac{\lambda^{2nk}}{4^{nk}},
\label{Eq:zeta-2m-n}\\
\zeta^{\star}(\{2k\}^n)=&\sum\limits_{n_0+\cdots+n_{k-1}=nk\atop n_j\geqslant 0}\frac{\beta_{2n_0}\cdots\beta_{2n_{k-1}}\rho_k^{\sum\limits_{j=0}^{k-1}2jn_j}}{(2n_0)!\cdots(2n_{k-1})!}\frac{\lambda^{2nk}}{4^{nk}},
\label{Eq:zeta-star-2m-n}
\end{align}
where $\rho_k=e^{\frac{\pi\sqrt{-1}}{k}}$, $\lambda=2\pi\sqrt{-1}$ and $\beta_{2m}=B_{2m}(2-4^m)$. Then an evaluation formula for $\zeta^t(\{2k\}^n)$ is given in the following theorem.

\begin{thm}
For any $k\in\mathbb{N}$ and any $n\in\mathbb{Z}_{\geqslant 0}$, we have
\begin{align}
\zeta^t(\{2k\}^n)=\sum\limits_{i+j=n\atop i,j\geqslant 0}\zeta(\{2k\}^i)\zeta^{\star}(\{2k\}^j)(1-t)^it^j
\label{Eq:zetat-2k-star}
\end{align}
and
\begin{align}
\zeta^t(\{2k\}^n)=&\sum\limits_{i+j=n\atop i,j\geqslant 0}(t-1)^{i}t^{j}\sum\limits_{{n_0+\cdots+n_{k-1}=ik\atop m_0+\cdots+m_{k-1}=jk}\atop {n_0,\ldots,n_{k-1}\geqslant 0\atop m_0,\ldots,m_{k-1}\geqslant 0}}\rho_k^{2\sum\limits_{l=0}^{k-1}l(n_l+m_l)}\nonumber\\
&\qquad\times\prod\limits_{l=0}^{k-1}
\frac{\beta_{2m_l}}{(2n_l+1)!(2m_l)!}\frac{\lambda^{2nk}}{4^{nk}},
\label{Eq:zetat-2k-power}
\end{align}
where $\lambda=2\pi\sqrt{-1}$, $\rho_k=e^{\frac{\pi\sqrt{-1}}{k}}$ and $\beta_{2m}=(2-4^m)B_{2m}$.
\end{thm}

\proof
Using Euler's formula
$$\zeta(2n)=-\frac{B_{2n}}{2(2n)!}\lambda^{2n},$$
which can be deduced from the extended double shuffle relations by \cite{Li-Qin}, we have
\begin{align*}
\sum\limits_{n=1}^\infty\frac{1}{n}d_n(t)\zeta(2nk)u^{2nk}=&\sum\limits_{n=1}^\infty\frac{B_{2nk}}{2n(2nk)!}(t-1)^nu^{2nk}\lambda^{2nk}\\
&-\sum\limits_{n=1}^\infty\frac{B_{2nk}}{2n(2nk)!}t^nu^{2nk}\lambda^{2nk}.
\end{align*}
Applying \eqref{Eq:Sum-BernoulliN}, we find
\begin{align*}
\sum\limits_{n=1}^\infty\frac{1}{n}d_n(t)\zeta(2nk)u^{2nk}=&\sum\limits_{j=0}^{k-1}\left(\log\frac{e^{u_1\rho_k^j}-1}{u_1\rho_k^j}-\frac{1}{2}u_1\rho_k^j\right)\\
&-\sum\limits_{j=0}^{k-1}\left(\log\frac{e^{u_2\rho_k^j}-1}{u_2\rho_k^j}-\frac{1}{2}u_2\rho_k^j\right),
\end{align*}
where $u_1=(t-1)^{\frac{1}{2k}}u\lambda$ and $u_2=t^{\frac{1}{2k}}u\lambda$. Hence we get
\begin{align}
\exp\left(\sum\limits_{n=1}^\infty\frac{1}{n}d_n(t)\zeta(2nk)u^{2nk}\right)=&\prod\limits_{j=0}^{k-1}\frac{e^{\frac{1}{2}u_1\rho_k^j}-e^{-\frac{1}{2}u_1\rho_k^j}}{u_1\rho_k^j}\nonumber\\
&\times\prod\limits_{j=0}^{k-1}\left(\frac{u_2\rho_k^j}{e^{\frac{1}{2}u_2\rho_k^j}-1}-\frac{u_2\rho_k^j}{e^{u_2\rho_k^j}-1}\right).
\label{Eq:zeta-2nk}
\end{align}

Using \cite[Proposition 3.17]{Li-Qin}, we have
$$\sum\limits_{n=0}^\infty\zeta(\{2k\}^n)(1-t)^nu^{2nk}=\prod\limits_{j=0}^{k-1}\frac{e^{\frac{1}{2}u_1\rho_k^j}-e^{-\frac{1}{2}u_1\rho_k^j}}{u_1\rho_k^j}$$
and
$$\sum\limits_{n=0}^\infty\zeta^{\star}(\{2k\}^n)t^nu^{2nk}=\prod\limits_{j=0}^{k-1}\left(\frac{u_2\rho_k^j}{e^{\frac{1}{2}u_2\rho_k^j}-1}-\frac{u_2\rho_k^j}{e^{u_2\rho_k^j}-1}\right).$$
Then \eqref{Eq:zetat-2k-star} follows from \eqref{Eq:zetat-zk-Power} and \eqref{Eq:zeta-2nk}. At last, one can get \eqref{Eq:zetat-2k-power} from \eqref{Eq:zeta-2m-n}, \eqref{Eq:zeta-star-2m-n} and \eqref{Eq:zetat-2k-star}.
\qed

Setting $t=0$ and $t=1$ in \eqref{Eq:zetat-2k-power}, we get \cite[Theorem 3.15]{Li-Qin}. Let $k=1$ and $k=2$ in \eqref{Eq:zetat-2k-power}, we get
$$\zeta^t(\{2\}^n)=\sum\limits_{i+j=n\atop i,j\geqslant 0}(t-1)^it^j
\frac{\beta_{2j}}{(2i+1)!(2j)!}\frac{\lambda^{2n}}{4^{n}}$$
and
$$\zeta^t(\{4\}^n)=\sum\limits_{i+j=n\atop i,j\geqslant 0}(t-1)^it^j\sum\limits_{i_1+i_2=2i,j_1+j_2=2j\atop i_1,i_2,j_1,j_2\geqslant 0}\frac{(-1)^{i_2+j_2}\beta_{2j_1}\beta_{2j_2}}{(2i_1+1)!(2i_2+1)!(2j_1)!(2j_2)!}\frac{\lambda^{4n}}{4^{2n}}.$$
Since
$$\sum\limits_{i=0}^{2n}(-1)^i\binom{4n+2}{2i+1}=(-1)^n2^{2n+1},$$
we find
$$\zeta^t(\{4\}^n)=\sum\limits_{i+j=n\atop i,j\geqslant 0}(1-t)^it^j\frac{2\lambda^{4n}}{(4i+2)!4^{n+j}}\sum\limits_{j_1+j_2=2j\atop j_1,j_2\geqslant 0}\frac{(-1)^{j_2}\beta_{2j_1}\beta_{2j_2}}{(2j_1)!(2j_2)!}.$$


\subsection{Some formulas of $\zeta^t(\{k\}^n)$}

We shall supply another method to prove \eqref{Eq:zetat-2k-star}. In fact, we can show the corresponding identity in the algebra $\mathfrak{h}_t^1$.

\begin{prop}\label{Prop:St-zt-power}
For any $k\in\mathbb{N}$, we have
\begin{align}
S_t\left(\frac{1}{1-z_ku}\right)=\frac{1}{1-z_k(1-t)u}\ast S\left(\frac{1}{1-z_ktu}\right),
\label{Eq:St-zku}
\end{align}
or equivalently
\begin{align}
S_t(z_k^n)=\sum\limits_{i+j=n\atop i,j\geqslant 0}(1-t)^it^jz_k^i\ast S(z_k^j),\quad (\forall n\geqslant 1).
\label{Eq:St-zk-power}
\end{align}
In particular, if $k\geqslant 2$, it holds
$$\zeta^t(\{k\}^n)=\sum\limits_{i+j=n\atop i,j\geqslant 0}(1-t)^it^j\zeta(\{k\}^i)\zeta^{\star}(\{k\}^j),\quad (\forall n\geqslant 1).$$
\end{prop}

To prove \eqref{Eq:St-zk-power}, we need the following lemma, which generalizes \cite[Eq. (4)]{Muneta2008}.

\begin{lem}
For any $n,k_1,\ldots,k_n\in\mathbb{N}$, we have
\begin{align}
S_t(z_{k_1}z_{k_2}\cdots z_{k_n})=\sum\limits_{i=1}^{n}t^{i-1}z_{k_1+\cdots+k_i}S_t(z_{k_{i+1}}\cdots z_{k_n}).
\label{Eq:St-Extension}
\end{align}
In particular, for any $k,n\in\mathbb{N}$, we have
\begin{align}
S_t(z_k^n)=\sum\limits_{i=1}^nt^{i-1}z_{ik}S_t(z_k^{n-i}).
\label{Eq:St-zk-induction}
\end{align}
\end{lem}

\proof
We prove  \eqref{Eq:St-Extension} by induction on. The case of $n=1$ is trivial. Now assume that $n\geqslant 2$. Since
$$S_t(z_{k_1}z_{k_2}\cdots z_{k_n})=z_{k_1}S_t(z_{k_2}\cdots z_{k_n})+tx^{k_1}S_t(z_{k_2}\cdots z_{k_n}),$$
using the inductive hypothesis, we get
$$S_t(z_{k_1}z_{k_2}\cdots z_{k_n})=z_{k_1}S_t(z_{k_2}\cdots z_{k_n})+tx^{k_1}\sum\limits_{i=2}^{n}t^{i-2}z_{k_2+\cdots+k_i}S_t(z_{k_{i+1}}\cdots z_{k_n}),$$
which finishes the proof.
\qed

Now we return to prove Proposition \ref{Prop:St-zt-power}.

\noindent {\bf Proof of Proposition \ref{Prop:St-zt-power}.}
We prove \eqref{Eq:St-zk-power} by induction on $n$. The case of $n=1$ is trivial. Now assume that $n\geqslant 2$. The right-hand side of \eqref{Eq:St-zk-power} is
$$(1-t)^nz_k^n+t^nS(z_k^n)+\sum\limits_{j=1}^{n-1}(1-t)^{n-j}t^jz_k^{n-j}\ast S(z_k^j).$$
Using \eqref{Eq:St-zk-induction} with $t=1$, we find the right-hand side of \eqref{Eq:St-zk-power} becomes
$$(1-t)^nz_k^n+t^nS(z_k^n)+\sum\limits_{j=1}^{n-1}\sum\limits_{i=1}^j(1-t)^{n-j}t^{j}z_k^{n-j}\ast z_{ik}S(z_k^{j-i}),$$
which is $t^nS(z_k^n)+T_1+T_2+T_3$, where
\begin{align*}
T_1=&(1-t)^nz_k^n+z_k\sum\limits_{j=1}^{n-1}\sum\limits_{i=1}^j(1-t)^{n-j}t^{j}z_k^{n-j-1}\ast z_{ik}S(z_k^{j-i}),\\
T_2=&\sum\limits_{j=1}^{n-1}\sum\limits_{i=1}^j(1-t)^{n-j}t^{j}z_{ik}\left[z_k^{n-j}\ast S(z_k^{j-i})\right],\\
T_3=&\sum\limits_{j=1}^{n-1}\sum\limits_{i=1}^j(1-t)^{n-j}t^{j}z_{(i+1)k}\left[z_k^{n-j-1}\ast S(z_k^{j-i})\right],
\end{align*}
respectively. By \eqref{Eq:St-zk-induction} with $t=1$, we have
$$T_1=(1-t)^nz_k^n+z_k\sum\limits_{j=1}^{n-1}(1-t)^{n-j}t^{j}z_k^{n-1-j}\ast S(z_k^j).$$
Then using the inductive hypothesis, we get
$$T_1=(1-t)z_kS_t(z_k^{n-1}).$$
Changing the order of the summation, we have
\begin{align*}
T_2=&\sum\limits_{i=1}^{n-1}t^{i}z_{ik}\sum\limits_{j=i}^{n-1}(1-t)^{n-j}t^{j-i}z_k^{n-j}\ast S(z_k^{j-i})\\
=&\sum\limits_{i=1}^{n-1}t^{i}z_{ik}\sum\limits_{j=0}^{n-i-1}(1-t)^{n-i-j}t^{j}z_k^{n-i-j}\ast S(z_k^{j}).
\end{align*}
While by the inductive hypothesis, we get
\begin{align*}
T_2=&\sum\limits_{i=1}^{n-1}t^{i}z_{ik}\left[S_t(z_k^{n-i})-t^{n-i}S(z_k^{n-i})\right]\\
=&\sum\limits_{i=1}^{n-1}t^{i}z_{ik}S_t(z_k^{n-i})-\sum\limits_{i=1}^{n-1}t^{n}z_{ik}S(z_k^{n-i}).
\end{align*}
Similarly, we have
\begin{align*}
T_3=&\sum\limits_{i=1}^{n-1}t^{i}z_{(i+1)k}\sum\limits_{j=i}^{n-1}(1-t)^{n-j}t^{j-i}z_k^{n-j-1}\ast S(z_k^{j-i})\\
=&\sum\limits_{i=1}^{n-1}t^{i}z_{(i+1)k}\sum\limits_{j=0}^{n-i-1}(1-t)^{n-i-j}t^{j}z_k^{n-i-1-j}\ast S(z_k^{j})\\
=&\sum\limits_{i=1}^{n-1}t^{i}z_{(i+1)k}(1-t)S_t(z_k^{n-i-1})=\sum\limits_{i=2}^{n}t^{i-1}z_{ik}(1-t)S_t(z_k^{n-i})\\
=&\sum\limits_{i=2}^{n}t^{i-1}z_{ik}S_t(z_k^{n-i})-\sum\limits_{i=2}^{n}t^{i}z_{ik}S_t(z_k^{n-i}).
\end{align*}
Hence the right-hand side of \eqref{Eq:St-zk-power} becomes
$$t^nS(z_k^n)-t^n\sum\limits_{i=1}^{n-1}z_{ik}S(z_k^{n-i})+\sum\limits_{i=1}^{n}t^{i-1}z_{ik}S_t(z_k^{n-i})-t^nz_{nk},$$
which is $S_t(z_k^n)$ by \eqref{Eq:St-zk-induction}.
\qed

In \cite[Corollary 1]{Ihara-Kajikawa-Ohno-Okuda}, it was proved that for any $k\in\mathbb{N}$, it holds
\begin{align}
S\left(\frac{1}{1-z_ku}\right)\ast\frac{1}{1+z_ku}=1.
\label{Eq:S-zp-ast-zp}
\end{align}
We have some generalizations of \eqref{Eq:S-zp-ast-zp}.

\begin{cor}
For any $k\in\mathbb{N}$, we have
\begin{align}
S_t\left(\frac{1}{1-z_ku}\right)\ast S_{1-t}\left(\frac{1}{1+z_ku}\right)=1
\label{Eq:St-S1-t}
\end{align}
and
\begin{align}
S_{1-2t}\left(\frac{1}{1-z_ku}\right)\tast \frac{1}{1+z_ku}=1.
\label{Eq:S1-2t-tast}
\end{align}
\end{cor}

\proof
Using \eqref{Eq:St-zku}, we get
$$S_{1-t}\left(\frac{1}{1+z_ku}\right)=\frac{1}{1+z_ktu}\ast S\left(\frac{1}{1+z_k(1-t)u}\right),$$
which together with \eqref{Eq:St-zku} and \eqref{Eq:S-zp-ast-zp} implies  \eqref{Eq:St-S1-t}. Applying $S_{-t}$ to \eqref{Eq:St-S1-t} and using Lemma \ref{Lem:S-Com}, we get \eqref{Eq:S1-2t-tast}.
\qed

Another possible generalization of \eqref{Eq:S-zp-ast-zp} is the following proposition.

\begin{prop}
For any $k\in\mathbb{N}$, we have
\begin{align}
S_t\left(\frac{1}{1-z_ku}\right)\ast\frac{1}{1+z_ku}=\frac{1}{1-\sum\limits_{i=2}^\infty t^{i-2}(t-1)z_{ik}u^i}.
\label{Eq:St-zp-ast-zp}
\end{align}
\end{prop}

\proof
Denote the left-hand side of \eqref{Eq:St-zp-ast-zp} by $A(u)$, then we get
$$A(u)=\sum\limits_{n=0}^\infty A_nu^n,$$
where
$$A_n=\sum\limits_{j=0}^n(-1)^jz_k^j\ast S_t(z_k^{n-j}).$$
For example, we have $A_0=1$, $A_1=0$ and $A_2=(t-1)z_{2k}$. Assume that $n\geqslant 2$. Then using \eqref{Eq:St-zk-induction}, we get
\begin{align*}
A_n=&S_t(z_k^n)+(-1)^nz_k^n+\sum\limits_{j=1}^{n-1}\sum\limits_{i=1}^{n-j}(-1)^jt^{i-1}z_k^j\ast z_{ik}S_t(z_k^{n-j-i})\\
=&S_t(z_k^n)+T_1+T_2+T_3,
\end{align*}
where
\begin{align*}
T_1=&(-1)^nz_k^n+z_k\sum\limits_{j=1}^{n-1}(-1)^jz_k^{j-1}\ast S_t(z_k^{n-j}),\\
T_2=&\sum\limits_{j=1}^{n-1}\sum\limits_{i=1}^{n-j}(-1)^jt^{i-1}z_{ik}\left(z_k^j\ast S_t(z_k^{n-j-i})\right),\\
T_3=&\sum\limits_{j=1}^{n-1}\sum\limits_{i=1}^{n-j}(-1)^jt^{i-1}z_{(i+1)k}\left(z_k^{j-1}\ast S_t(z_k^{n-j-i})\right),
\end{align*}
respectively. Obviously, one has
$$T_1=z_k\sum\limits_{j=1}^{n}(-1)^jz_k^{j-1}\ast S_t(z_k^{n-j})=-z_kA_{n-1}.$$
Changing the order of the summation, we have
\begin{align*}
T_2=&\sum\limits_{i=1}^{n-1}t^{i-1}z_{ik}\sum\limits_{j=1}^{n-i}(-1)^{j}z_k^j\ast S_t(z_k^{n-i-j})\\
=&\sum\limits_{i=1}^{n-1}t^{i-1}z_{ik}\left(A_{n-i}-S_t(z_k^{n-i})\right).
\end{align*}
Hence by \eqref{Eq:St-zk-induction}, we get
$$T_2=\sum\limits_{i=1}^{n-1}t^{i-1}z_{ik}A_{n-i}-S_t(z_k^n)+t^{n-1}z_{nk}.$$
Similarly, we have
\begin{align*}
T_3=&\sum\limits_{i=1}^{n-1}t^{i-1}z_{(i+1)k}\sum\limits_{j=1}^{n-i}(-1)^jz_k^{j-1}\ast S_t(z_k^{n-i-j})\\
=&-\sum\limits_{i=1}^{n-1}t^{i-1}z_{(i+1)k}A_{n-i-1}=-\sum\limits_{i=2}^{n}t^{i-2}z_{ik}A_{n-i}.
\end{align*}
Therefore we obtain
$$A_n=\sum\limits_{i=2}^n(t^{i-1}-t^{i-2})z_{ik}A_{n-i}.$$
Applying $\sum\limits_{n=2}^\infty u^n$ to both sides of the above equality,  we get
\begin{align*}
A(u)-1=&\sum\limits_{n=2}^\infty\sum\limits_{i=2}^n(t^{i-1}-t^{i-2})z_{ik}A_{n-i}u^n\\
=&\sum\limits_{i=2}^\infty(t^{i-1}-t^{i-2})z_{ik}u^i\sum\limits_{n=i}^\infty A_{n-i}u^{n-i}\\
=&\sum\limits_{i=2}^\infty(t^{i-1}-t^{i-2})z_{ik}u^iA(u),
\end{align*}
which implies the desired result.
\qed

Setting $t=1$ in \eqref{Eq:St-zp-ast-zp}, we get \eqref{Eq:S-zp-ast-zp}. Setting $t=0$ in \eqref{Eq:St-zp-ast-zp}, we find
$$\frac{1}{1-z_ku}\ast\frac{1}{1+z_ku}=\frac{1}{1+z_{2k}u^2},$$
which is equivalent to
$$\sum\limits_{j=0}^n (-1)^jz_k^{j}\ast z_k^{n-j}=\begin{cases}
0 & \text{if\;} n\text{\, is odd},\\
(-1)^mz_{2k}^m & \text{if\;} n=2m \text{\, is even}.
\end{cases}$$
The above equality for odd integer $n$  is  trivial, while for even integer $n$ might be interesting. In particular, for any integers $k,n$ with $k\geqslant 2$ and $n\geqslant 0$, we have
$$\zeta(\{2k\}^n)=\sum\limits_{j=0}^{2n}(-1)^{n-j}\zeta(\{k\}^j)\zeta(\{k\}^{2n-j}),$$
or equivalently
$$\zeta(\{2k\}^n)=2\sum\limits_{j=0}^{n-1}(-1)^{n-j}\zeta(\{k\}^j)\zeta(\{k\}^{2n-j})+\zeta(\{k\}^n)^2.$$

Using \eqref{Eq:S-zp-ast-zp} and \eqref{Eq:St-zp-ast-zp}, we get the following corollary immediately.

\begin{cor}
Let $k\in\mathbb{N}$. We have
$$S_t\left(\frac{1}{1-z_ku}\right)=\frac{1}{1-\sum\limits_{i=2}^\infty t^{i-2}(t-1)z_{ik}u^i}\ast S\left(\frac{1}{1-z_ku}\right).$$
In particular, for any $n\in\mathbb{Z}_{\geqslant 0}$, it holds
\begin{align}
S_t(z_k^n)=S(z_k^n)+\sum\limits_{j=2}^n\sum\limits_{m=1}^{\left[\frac{j}{2}\right]}t^{j-2m}(t-1)^m\left(\sum\limits_{i_1+\cdots+i_m=j\atop i_1,\ldots,i_m\geqslant 2}z_{i_1k}\cdots z_{i_mk}\right)\ast S(z_k^{n-j}).
\label{Eq:St-zk-power-new}
\end{align}
And if $k\geqslant 2$, it holds
\begin{align}
\zeta^t(\{k\}^n)=&\zeta^{\star}(\{k\}^n)\nonumber\\
&+\sum\limits_{j=2}^n\sum\limits_{m=1}^{\left[\frac{j}{2}\right]}t^{j-2m}(t-1)^m\left(\sum\limits_{i_1+\cdots+i_m=j\atop i_1,\ldots,i_m\geqslant 2}\zeta(i_1k,\ldots, i_mk)\right)\zeta^{\star}(\{k\}^{n-j}).
\label{Eq:zetat-k-zeta-star}
\end{align}
\end{cor}

As an application of the above corollary, setting $t=0$ in \eqref{Eq:zetat-k-zeta-star}, we obtain
$$\zeta(\{k\}^n)=\sum\limits_{m=0}^{\left[\frac{n}{2}\right]}(-1)^m\zeta(\{2k\}^m)\zeta^{\star}(\{k\}^{n-2m}),\quad (k\geqslant 2,n\geqslant 0).$$


\subsection{Restricted sum formula}

In this subsection, let $a$ be a fixed positive integer. For any $k,n\in\mathbb{N}$ with $k\geqslant n$, set
$$N_{k,n}=\sum\limits_{k_1+\cdots+k_n=k\atop k_i\geqslant 1}z_{ak_1}\cdots z_{ak_n}\in\mathfrak{h}^1.$$
We are interested in $S_t(N_{k,n})$. One can represent $S_t(N_{k,n})$ by $N_{k,i}$ for $i=1,\ldots,n$, which generalizes \cite[Eq. (3.41)]{Li-Qin}.

\begin{prop}
For any $k,n\in\mathbb{N}$ with $k\geqslant n$, we have
\begin{align}
S_t(N_{k,n})=\sum\limits_{i=1}^n\binom{k-i}{k-n}t^{n-i}N_{k,i}.
\label{Eq:St-Nkn-First}
\end{align}
\end{prop}

\proof
We prove \eqref{Eq:St-Nkn-First} by induction on $n$. The case of $n=1$ is trivial. Now assume that $k\geqslant n>1$. Then by the definition of $S_t$, we have
\begin{align*}
S_t(N_{k,n})=&\sum\limits_{k_1+\cdots+k_n=k\atop k_i\geqslant 1}(x^{ak_1-1}y+tx^{ak_1})S_t(z_{ak_2}\cdots z_{ak_n})\\
=&\sum\limits_{l=1}^{k-n+1}(x^{al-1}y+tx^{al})S_t(N_{k-l,n-1}).
\end{align*}
By the inductive assumption, we get
\begin{align*}
S_t(N_{k,n})=&\sum\limits_{l=1}^{k-n+1}\sum\limits_{i=1}^{n-1}\binom{k-l-i}{k-l-n+1}t^{n-1-i}z_{al}N_{k-l,i}\\
&+\sum\limits_{l=1}^{k-n+1}\sum\limits_{i=1}^{n-1}\binom{k-l-i}{k-l-n+1}t^{n-i}x^{al}N_{k-l,i}.
\end{align*}
Since the first term of the right-hand side of the above equation is
$$\sum\limits_{i=2}^n\sum\limits_{k_1+\cdots+k_i=k\atop k_j\geqslant 1}\binom{k-k_1-i+1}{n-i}t^{n-i}z_{ak_1}\cdots z_{ak_i},$$
and the second term is
$$\sum\limits_{i=1}^{n-1}\sum\limits_{k_1+\cdots+k_i=k\atop k_j\geqslant 1}\left(\sum\limits_{l=1}^{k_1-1}\binom{k-l-i}{n-i-1}\right)t^{n-i}z_{ak_1}\cdots z_{ak_i},$$
we get
\begin{align*}
S_t(N_{k,n})=&\sum\limits_{k_1+\cdots+k_n=k\atop k_j\geqslant 1}z_{ak_1}\cdots z_{ak_n}+\sum\limits_{l=1}^{k-1}\binom{k-l-1}{n-2}t^{n-1}z_{ak}\\
&+\sum\limits_{i=2}^{n-1}\sum\limits_{k_1+\cdots+k_i=k\atop k_j\geqslant 1}\left(\sum\limits_{l=1}^{k_1-1}\binom{k-l-i}{n-i-1}+\binom{k-k_1-i+1}{n-i}\right)t^{n-i}z_{ak_1}\cdots z_{ak_i}\\
=&\sum\limits_{k_1+\cdots+k_n=k\atop k_j\geqslant 1}z_{ak_1}\cdots z_{ak_n}+\binom{k-1}{n-1}t^{n-1}z_{ak}\\
&+\sum\limits_{i=2}^{n-1}\sum\limits_{k_1+\cdots+k_i=k\atop k_j\geqslant 1}\binom{k-i}{n-i}t^{n-i}z_{ak_1}\cdots z_{ak_i},
\end{align*}
as required.
\qed

Since
$$N_{k,i}=\sum\limits_{j=0}^{k-i}(-1)^{k-i-j}\binom{k-j}{i}S(z_a^j)\ast z_a^{k-j},$$
which was proved in \cite[Proposition 3.26]{Li-Qin}, we get the following result.

\begin{cor}
For any $k,n\in\mathbb{N}$ with $k\geqslant n$, we have
\begin{align}
S_t(N_{k,n})=\sum\limits_{j=0}^{k-1}\left[\sum\limits_{i=1}^{\min\{n,k-j\}}(-1)^{k-i-j}\binom{k-i}{k-n}\binom{k-j}{i}t^{n-i}\right]S(z_a^j)\ast z_a^{k-j}.
\label{Eq:St-Nkn-Second}
\end{align}
In particular, if $a>1$, we have
\begin{align}
\sum\limits_{k_1+\cdots+k_n=k\atop k_i\geqslant 1}\zeta^t(ak_1,\ldots,ak_n)=&\sum\limits_{j=0}^{k-1}\left[\sum\limits_{i=1}^{\min\{n,k-j\}}(-1)^{k-i-j}\binom{k-i}{k-n}\binom{k-j}{i}t^{n-i}\right]\nonumber\\
&\quad\times\zeta^{\star}(\{a\}^j)\zeta(\{a\}^{k-j}).
\label{Eq:zetat-restricted-sum}
\end{align}
\end{cor}

Then using the evaluation formulas \eqref{Eq:zeta-2m-n} and \eqref{Eq:zeta-star-2m-n}, we get the following restricted sum formula.

\begin{thm}
For any $k,n,m\in\mathbb{N}$ with $k\geqslant n$, we have
\begin{align*}
&\sum\limits_{k_1+\cdots+k_n=k\atop k_i\geqslant 1}\zeta^t(2mk_1,\ldots,2mk_n)\\
=&\left\{\sum\limits_{i=1}^{n}\left[\sum\limits_{j=0}^{k-i}\binom{k-j}{i}\sum\limits_{{n_0+\cdots+n_{m-1}=mj\atop l_0+\cdots+l_{m-1}=m(k-j)}\atop n_p,l_p\geqslant 0}\prod\limits_{p=0}^{m-1}\frac{\beta_{2n_p}}{(2n_p)!(2l_p+1)!}\right.\right.\\
&\qquad\left.\left.\times \rho_m^{\sum\limits_{p=0}^{m-1}2p(n_p+l_p)}\right](-1)^i\binom{k-i}{k-n}t^{n-i}\right\}\frac{\lambda^{2km}}{4^{km}},
\end{align*}
where $\rho_m=e^{\frac{\pi\sqrt{-1}}{m}}$, $\lambda=2\pi\sqrt{-1}$ and $\beta_{2l}=B_{2l}(2-4^l)$.
\end{thm}

\proof
Using \eqref{Eq:zeta-2m-n}, \eqref{Eq:zeta-star-2m-n} and \eqref{Eq:zetat-restricted-sum}, we get
\begin{align*}
&\sum\limits_{k_1+\cdots+k_n=k\atop k_i\geqslant 1}\zeta^t(2mk_1,\ldots,2mk_n)=\sum\limits_{j=0}^{k-1}\left[\sum\limits_{i=1}^{\min\{n,k-j\}}(-1)^i\binom{k-i}{k-n}\binom{k-j}{i}t^{n-i}\right]\\
&\qquad\times\left[\sum\limits_{{n_0+\cdots+n_{m-1}=mj\atop l_0+\cdots+l_{m-1}=m(k-j)}\atop n_p,l_p\geqslant 0}\prod\limits_{p=0}^{m-1}\frac{\beta_{2n_p}}{(2n_p)!(2l_p+1)!}\rho_m^{\sum\limits_{p=0}^{m-1}2p(n_p+l_p)}\right]\frac{\lambda^{2km}}{4^{km}}.
\end{align*}
Changing the order of the summation, we obtain the desired result.
\qed

Another way to get \eqref{Eq:St-Nkn-Second} is using the generating functions, which is similar as in \cite[Subsection 6.1]{Hoffman-Ihara}. We set
\begin{align*}
E(u)=&\sum\limits_{j=0}^\infty z_a^ju^j=\frac{1}{1-z_au},\\
H(u)=&\sum\limits_{j=0}^\infty S(z_a^j)u^j=S(E(u)),\\
\mathcal{F}(u,v)=&1+\sum\limits_{k\geqslant n\geqslant 1}N_{k,n}u^kv^n.
\end{align*}
Using \cite[Eq.(3.39)]{Li-Qin}, we have
$$\mathcal{F}(u,v)=E((v-1)u)\ast H(u).$$
We can rewrite \eqref{Eq:St-zku} as
$$S_t(E(u))=E((1-t)u)\ast H(tu).$$
Hence we get
$$\mathcal{F}(u,v)=S_{\frac{1}{v}}(E(uv)).$$
Using Lemma \ref{Lem:S-Com}, we have
$$S_t(\mathcal{F}(u,v))=S_{t+\frac{1}{v}}(E(uv)).$$
Then applying \eqref{Eq:St-zku} again, we obtain the following result.

\begin{prop}
We have
\begin{align}
S_t(\mathcal{F}(u,v))=E((v-tv-1)u)\ast H((1+tv)u).
\label{Eq:St-F}
\end{align}
\end{prop}

Comparing the coefficients of $u^kv^n$ of both sides of \eqref{Eq:St-F}, we easily get another expression of $S_t(N_{k,n})$.

\begin{cor}
For any $k,n\in\mathbb{N}$ with $k\geqslant n$, we have
\begin{align}
S_t(N_{k,n})=\sum\limits_{j_1+j_2=k\atop j_1,j_2\geqslant 0}\left[\sum\limits_{{i_1+i_2=n\atop 0\leqslant i_1\leqslant j_1}\atop 0\leqslant i_2\leqslant j_2}(-1)^{j_1-i_1}\binom{j_1}{i_1}\binom{j_2}{i_2}(1-t)^{i_1}t^{i_2}\right]z_a^{j_1}\ast S(z_a^{j_2}).
\label{Eq:St-Nkn-Third}
\end{align}
\end{cor}

Finally, we show that \eqref{Eq:St-Nkn-Third} is equivalent to \eqref{Eq:St-Nkn-Second}. In fact, the inner sum of the right-hand side of \eqref{Eq:St-Nkn-Third} is
$$\sum\limits_{{i_1+i_2=n\atop 0\leqslant l\leqslant i_1\leqslant j_1}\atop 0\leqslant i_2\leqslant j_2}(-1)^{j_1-i_1+l}\binom{j_1}{i_1}\binom{j_2}{i_2}\binom{i_1}{l}t^{l+i_2},$$
which equals
\begin{align*}
&\sum\limits_{{i_1+i_2=n\atop 0\leqslant i\leqslant i_1\leqslant j_1}\atop 0\leqslant i_2\leqslant j_2}(-1)^{j_1+i}\binom{j_1}{i_1}\binom{j_2}{i_2}\binom{i_1}{i}t^{n-i}\\
=&\sum\limits_{{0\leqslant i\leqslant i_1\leqslant j_1}\atop 0\leqslant n-i_1\leqslant j_2}(-1)^{j_1+i}\binom{j_1}{i}\binom{j_1-i}{i_1-i}\binom{j_2}{n-i_1}t^{n-i}\\
=&\sum\limits_{i=0}^{\min\{n,j_1\}}(-1)^{j_1+i}\binom{j_1}{i}t^{n-i}\sum\limits_{i\leqslant i_1\leqslant j_1\atop n-j_2\leqslant i_1\leqslant n}\binom{j_1-i}{i_1-i}\binom{j_2}{n-i_1}.
\end{align*}
Using the Chu-Vandermonde identity, we find the inner sum of the right-hand side of \eqref{Eq:St-Nkn-Third} is
$$\sum\limits_{i=0}^{\min\{n,j_1\}}(-1)^{j_1+i}\binom{j_1}{i}\binom{j_1+j_2-i}{n-i}t^{n-i}.$$
Hence the right-hand side of \eqref{Eq:St-Nkn-Third} is
$$\sum\limits_{j=0}^k\left[\sum\limits_{i=0}^{\min\{n,k-j\}}(-1)^{k-j-i}\binom{k-j}{i}\binom{k-i}{n-i}t^{n-i}\right]z_a^{k-j}\ast S(z_a^j),$$
which equals
\begin{align*}
&\sum\limits_{j=0}^k(-1)^{k-j}\binom{k}{n}t^nz_a^{k-j}\ast S(z_a^j)\\
&+\sum\limits_{j=0}^{k-1}\left[\sum\limits_{i=1}^{\min\{n,k-j\}}(-1)^{k-j-i}\binom{k-j}{i}\binom{k-i}{k-n}t^{n-i}\right]z_a^{k-j}\ast S(z_a^j).
\end{align*}
While \eqref{Eq:S-zp-ast-zp} implies that
$$\sum\limits_{j=0}^k(-1)^{k-j}z_a^{k-j}\ast S(z_a^j)=0$$
for any $k\in\mathbb{N}$. Hence the right-hand side of \eqref{Eq:St-Nkn-Third} is exactly the same as the right-hand side of \eqref{Eq:St-Nkn-Second}.


\end{document}